\documentclass[a4paper,12pt]{article}
\usepackage{dsfont}
\usepackage{amsmath}
\usepackage{enumitem}
\usepackage{amssymb,amsthm}
\usepackage{mathrsfs,amsfonts}
\usepackage{graphicx}
\usepackage{bbm}
\usepackage{ulem}
\usepackage[hidelinks]{hyperref}
\usepackage[capitalise]{cleveref}

\allowdisplaybreaks

\newtheorem{thm0}{Theorem}[section]
\newtheorem{lem1}[thm0]{Lemma}
\newtheorem{thm1}[thm0]{Theorem}
\newtheorem{cor1}[thm0]{Corollary}
\newtheorem{pro1}[thm0]{Proposition}
\newtheorem{rem1}[thm0]{Remark}
\newtheorem{crit1}[thm0]{Criterion}

\def\bglemma{\begin{lem1}}\def\edlemma{\end{lem1}}
\def\bgtheorem{\begin{thm1}}\def\edtheorem{\end{thm1}}
\def\bgproposition{\begin{pro1}}\def\edproposition{\end{pro1}}
\def\bgremark{\begin{rem1}\rm{}\def\edremark{\end{rem1}}}

\def\beqlb{\begin{eqnarray}}\def\eeqlb{\end{eqnarray}}
\def\beqnn{\begin{eqnarray*}}\def\eeqnn{\end{eqnarray*}}
\def\bitemize{\begin{itemize}}\def\benumerate{\begin{itemize}}
\def\eitemize{\end{itemize}}
\def\ar{\!\!\!&}
\def\proof{\noindent{\it Proof.~}}\def\qed{\hfill$\square$\smallskip}
\def\mrm{\mathrm}\def\mbb{\mathbb}
\def\d{\mrm{d}}\def\e{\mrm{e}}

\begin{document}
	
\begin{center}
{\large\bfseries
\linespread{1.15}\selectfont
Domains of Attraction for Quasi-Stationary\\
Distributions of Branching Processes\\
in Random Environments}
\end{center}

\bigskip

\centerline{Pei-Sen Li}
\smallskip
\centerline{\it School of Mathematics and Statistics, Beijing Institute of Technology}
\centerline{\it Beijing 100872, China}
\centerline{\tt peisenli@bit.edu.cn}

\bigskip

{\narrower{\narrower
\noindent{\textit{Abstract:}} We obtain a complete characterization of the
quasi-stationary distributions of subcritical stable continuous-state branching
processes in Brownian random environments, together with the domain of attraction
of each member. These distributions form a one-parameter family whose modified
Laplace transforms are given explicitly in terms of Tricomi's confluent
hypergeometric function. For every initial distribution, we determine whether
the {conditional law} converges to a QSD and identify that QSD whenever it does.
The necessary and sufficient criteria are formulated in terms of the tail
probabilities and moments of the initial law.\bigskip


\noindent{\textit{2020 Mathematics Subject Classification:}\newline
Primary~60J80; Secondary~60K37,~60J25.}
\par}\par}

\section{Introduction and main results}
\setcounter{equation}{0}

\subsection{Background}\label{sec_intro_background}

Consider a time-homogeneous non-negative Markov process
$X=(X_t)_{t\ge0}$ on $\mbb R_+=[0,\infty)$, with $0$ absorbing. Write
$\tau_0:=\inf\{t\ge0:X_t=0\}$ for its extinction time, with the convention
$\inf\emptyset=\infty$. In the subcritical setting considered below,
$\tau_0<\infty$ almost surely. Since {the unconditional law of $X_t$
converges weakly to the point mass at $0$ as $t\to\infty$}, we are interested
in the process conditioned on survival for a long time.

For $x\in\mbb R_+$, let $\mbb P_x$ and $\mbb E_x$ denote the law and
expectation of the process started from $x$. If $\mu$ is a probability law on
$(0,\infty)$, we write $\mbb P_\mu$ and $\mbb E_\mu$ for
{the law and expectation of the process with initial distribution $\mu$}.
We say that a probability measure $\mu$ on $(0,\infty)$ is a
quasi-stationary distribution (QSD) if, for every $t\ge0$ and every measurable
set $A\subset(0,\infty)$,
\[
    \mu(A)=\mbb P_\mu(X_t\in A\mid t<\tau_0).
\]
Thus a QSD describes a non-trivial distribution of the population size that
remains unchanged when the process is conditioned on non-extinction.

Define the killed sub-Markov
semigroup $(P_t^0)_{t\ge0}$ on $(0,\infty)$ by
\begin{equation}\label{def_killed_semigroup_intro}
    P_t^0f(x)
    :=
    \mbb E_x\left[f(X_t)\mathbf 1_{\{t<\tau_0\}}\right],
    \qquad x>0,
\end{equation}
for bounded measurable functions $f$ on $(0,\infty)$.  For a probability
measure $\mu$ on $(0,\infty)$, \eqref{def_killed_semigroup_intro} gives
\[
    \mu P_t^0(A)
    :=\int_{(0,\infty)}P_t^0\mathbf 1_A(x)\,\mu(\d x)
    =\mbb P_\mu(X_t\in A,\ t<\tau_0),
\]
for every measurable set $A\subset(0,\infty)$.  A probability measure
$\mu$ is a QSD if and only if it is a
positive left eigenmeasure of this
semigroup: there exists $\rho>0$ such that
\beqlb\label{qsd_left_eigenmeasure}
    \mu P_t^0=\e^{-\rho t}\mu,
    \qquad t\ge0 .
\eeqlb
Consequently,
\[
    \mbb P_\mu(\tau_0>t)=\e^{-\rho t},
    \qquad t\geq0.
\]
Thus the extinction time $\tau_0$ under $\mbb P_\mu$ is exponentially
distributed with parameter $\rho$.
The QSD problem is therefore a Perron--Frobenius-type problem for the killed
semigroup: one seeks normalized positive left eigenmeasures and their
associated decay parameters.  In finite state spaces, the classical
Perron--Frobenius theorem gives, under standard assumptions, a unique
normalized positive left eigenvector and hence a unique QSD.  For processes
with non-compact state spaces, however, neither the existence nor the
uniqueness of a QSD is guaranteed.  The spectral problem is then to
characterize all QSDs and their decay parameters.  These distributions are
natural candidates for the long-time limits of the conditioned process.  One
is therefore led to ask whether such limits exist and how they depend on the
initial distribution.  Throughout, $\Longrightarrow$ denotes weak convergence
of probability measures.  For a fixed initial state $x>0$, if
\[
    \mbb P_x(X_t\in\cdot\mid t<\tau_0)
    \Longrightarrow \nu
    \qquad\text{as }t\to\infty,
\]
then $\nu$ is called a Yaglom limit from $x$.  When several QSDs are present,
different initial laws may select different members of the QSD family.
{The domain of attraction of a QSD consists of all initial laws $\mu$
for which the conditional law from $\mu$ converges weakly to that QSD.}

Quasi-stationarity in branching models goes back to Yaglom's pioneering limit
theorem for a subcritical Galton--Watson process \cite{Yag47}.  Seneta and
Vere-Jones \cite{SVJ66} subsequently identified a one-parameter family of QSDs
for this process, and Rubin and Vere-Jones \cite{RVJ68} determined the domains
of attraction of this family in terms of regular variation of the generating
function of the initial law.  A complete description of all QSDs of the
subcritical Galton--Watson process was later obtained by Maillard
\cite{Mai18}.  For
classical continuous-state branching processes, Lambert \cite{Lam07}
characterized the full family of QSDs under natural extinction assumptions and
showed that the Yaglom distribution is the stochastically smallest QSD in this
family.   The addition of competition can change this
picture.  When the competition mechanism is sufficiently strong near infinity,
the process may admit a unique QSD attracting all initial distributions; see
\cite{CCLMSM09,LWZ24,LZZ26SF} {for branching processes with competition}.
For birth--death processes absorbed almost surely,
van Doorn \cite{vD91} characterized the QSDs and determined the conditional
limits for finitely supported initial distributions. Conditions under which a
unique QSD attracts all initial distributions were also obtained for birth--death processes by
Zhang and Zhu \cite{ZZ13} and for one-dimensional diffusions by Zhang
and He \cite{ZH16}. For Brownian motion with negative drift absorbed at
zero, Mart{\'i}nez, Picco and San Mart{\'i}n
{\cite{MPSM98}} give sufficient conditions for attraction
in terms of the exponential decay of the initial density or tail; see also
Ben-Ari and Lee {\cite{BAL22}}.

The theory of
branching processes in random
environments and the extinction asymptotics relevant here were developed in
\cite{BH12,PP17,PPS16,LX18}.  {Li, Zheng and Zhou \cite{LZZ}
identified the common Yaglom limit for every fixed positive initial state
in each of the three subcritical regimes. This result leaves open the
classification of all QSDs and the characterization of their domains of
attraction for general initial laws.}

In this paper, we answer the two questions above for the stable
continuous-state branching process in a Brownian random environment (CBBRE).
First, we characterize all its QSDs and their decay parameters.  They form a
one-parameter family whose modified Laplace transforms are expressed in terms
of Kummer functions.  Second, we give a complete characterization of the attraction
problem for arbitrary initial laws.  {More precisely, every initial law
either belongs to the domain of attraction of exactly one QSD
or has a conditional law that converges to no QSD.
Explicit tail or moment criteria determine whether such convergence occurs and,
if so, which QSD is selected.}

\subsection{Model and main results}\label{sec_intro_main_results}

Let $(\Omega,\mathcal F,(\mathcal F_t)_{t\ge0},\mbb P)$ be a filtered
probability space.  We work with the stable continuous-state branching process
in a Brownian random environment whose branching mechanism is
\[
    \psi(\lambda)=-\alpha\lambda+c\lambda^{1+\beta},
    \qquad
    \lambda\ge0,
\]
where $c>0$, $\beta\in(0,1]$, and $\alpha\in\mbb R$.
Here $1+\beta$ is the stability index of the branching noise.
Equivalently, following the standard constructions in
\cite{BH12,PP17}, $Z=(Z_t)_{t\ge0}$ is the unique non-negative strong solution
of one of the following two stochastic equations, according to the value of
the stable index.  When $\beta=1$,
\beqlb\label{main_sde_feller}
Z_t
=
Z_0
+\alpha\int_0^t Z_s\,\d s
+\sigma\int_0^t Z_s\,\d B^{(e)}_s
+\sqrt{2c}\int_0^t\sqrt{Z_s}\,\d B^{(b)}_s ,
\eeqlb
where $B^{(e)}=(B_t^{(e)})_{t\geq0}$ and
$B^{(b)}=(B_t^{(b)})_{t\geq0}$ are independent standard Brownian motions and
$\sigma>0$ is the environmental volatility.  When $0<\beta<1$,
\beqlb\label{main_sde}
Z_t
=
Z_0
+\alpha\int_0^t Z_s\,\d s
+\sigma\int_0^t Z_s\,\d B^{(e)}_s
+\int_0^t\int_0^\infty\int_0^{Z_{s-}}
z\,\tilde N(\d s,\d z,\d u),
\eeqlb
where $(B_t^{(e)})_{t\geq0}$ is a standard Brownian motion and $\tilde N$ is an independent
compensated Poisson random measure with intensity
\[
    \frac{c\beta(\beta+1)}{\Gamma(1-\beta)}
    z^{-2-\beta}\,\d s\,\d z\,\d u .
\]

For $z\geq0$, let $\mbb P_z$ and $\mbb E_z$ denote the law and expectation
of $Z$ started from $z$. For a probability measure $\mu$ on $(0,\infty)$,
define
\[
    \mbb P_\mu(\cdot)
    :=\int_{(0,\infty)}\mbb P_z(\cdot)\,\mu(\d z),
\]
and let $\mbb E_\mu$ denote the corresponding expectation.
The state $0$ is absorbing.  We study the killed process on $(0,\infty)$ and
{refer to $\mbb P_\mu(Z_t\in\cdot\mid Z_t>0)$ as the conditional law
from $\mu$ at time $t$}.  Its killed
semigroup is again denoted by $(P_t^0)_{t\geq0}$.

The effective growth drift is defined by ${\bf m}:=\alpha-\sigma^2/2$.
Throughout the paper we assume the subcritical condition ${\bf m}<0$, so that
extinction occurs almost surely; see \cite[Theorem~5, arXiv version~2]{BH12} for $\beta=1$
and \cite[Theorem~3]{PP17} for $0<\beta<1$.  Set
\beqlb\label{def_kappa_c_front}
\kappa_c=-\frac{{\bf m}}{\sigma^2}>0.
\eeqlb
The survival asymptotics of subcritical branching processes in random
environments exhibit three regimes.  In the Brownian continuous-state
setting, this trichotomy was established for branching diffusions in
\cite{BH12} and for stable CBBREs in \cite{PP17}:
\[
\begin{array}{lll}
    0<\kappa_c<1
    &\quad&
    \text{weakly subcritical},\\[1mm]
    \kappa_c=1
    &&
    \text{intermediately subcritical},\\[1mm]
    \kappa_c>1
    &&
    \text{strongly subcritical}.
\end{array}
\]

Set $\kappa_*:=\min(1,\kappa_c)$. This is the upper endpoint of the
QSD parametrization and the critical moment order in the Yaglom
domain-of-attraction criterion below. For $\kappa\geq0$, define

\beqlb\label{gamma_kappa_front}
\gamma(\kappa)
=
-{\bf m}\kappa-\frac12\sigma^2\kappa^2.
\eeqlb
The function $\gamma$ is strictly increasing on the admissible interval
$(0,\kappa_*]$.  As a quadratic function on $[0,\infty)$, it attains its
maximum at $\kappa=\kappa_c$, where
$\gamma(\kappa_c)={\bf m}^2/(2\sigma^2)$.

For a probability measure $\mu$ on $[0,\infty)$, define its modified
Laplace transform by
\begin{equation}\label{def_modified_laplace_general}
    H_\mu(\lambda)=\int_0^\infty(1-\e^{-\lambda z})\,\mu(\d z),
    \qquad \lambda\geq0.
\end{equation}
To state the modified Laplace transforms of the QSDs explicitly, let $U(a,b,z)$ denote
Tricomi's confluent hypergeometric function. For $a>0$, $b\in\mathbb R$,
and $z>0$, it is given by
\beqlb\label{def_tricomi_U_integral}
    U(a,b,z)
    =\frac{1}{\Gamma(a)}
      \int_0^\infty \e^{-zr}r^{a-1}(1+r)^{b-a-1}\,\d r;
\eeqlb
see \cite[(13.4.4)]{NISTDLMF}. It solves Kummer's equation
$zw''+(b-z)w'-aw=0$ and satisfies $U(a,b,z)\sim z^{-a}$ as $z\to\infty$; see
\cite[(13.2.1), (13.2.6), and Section~13.7]{NISTDLMF}.
For $0<\kappa\leq\kappa_*$ and $\lambda>0$, set
\beqlb\label{def_Phi_kappa}
    \Phi_\kappa(\lambda)
    =
    \left(\frac{2c}{\beta\sigma^2}\lambda^\beta\right)^{\kappa/\beta}
    U\left(
        \frac{\kappa}{\beta},
        1+\frac{2(\kappa-\kappa_c)}{\beta},
        \frac{2c}{\beta\sigma^2}\lambda^\beta
    \right).
\eeqlb
Since $U(a,b,z)\sim z^{-a}$ as $z\to\infty$,
\[
    \lim_{\lambda\to\infty}\Phi_\kappa(\lambda)=1,
    \qquad 0<\kappa\leq\kappa_*.
\]

Our first theorem characterizes all QSDs of the CBBRE.

\bgtheorem[Complete QSD characterization]\label{thm_complete_qsd_classification}
The following two assertions hold.
\begin{itemize}
    \item[\rm(i)] If \(\nu\) is any QSD for the CBBRE with
    decay rate \(\rho>0\), then there exists a unique
    \(\kappa\in(0,\kappa_*]\) such that
    \[
        \rho=\gamma(\kappa),\qquad
        H_\nu(\lambda)=\Phi_\kappa(\lambda),\qquad \lambda>0.
    \]
    In particular, there is at most one QSD with each decay rate.
    \item[\rm(ii)] Conversely, for every \(0<\kappa\leq\kappa_*\),
    there exists a unique QSD with decay rate \(\gamma(\kappa)\),
    denoted by \(\nu_\kappa\).
\end{itemize}
\edtheorem

The member $\nu_{\kappa_*}$ is the Yaglom distribution identified in
\cite[Theorem~1.1]{LZZ}. More precisely, for every initial state $z>0$,
\[
    \mbb P_z(Z_t\in\cdot\mid Z_t>0)
    \Longrightarrow\nu_{\kappa_*}
    \qquad\text{as }t\to\infty.
\]
{We now turn from fixed initial states to general initial distributions
and study the corresponding conditional limits. The following theorem gives
a complete characterization of the domain of attraction of each QSD.}

Let $\mathcal P((0,\infty))$ denote the set of all probability measures on
$(0,\infty)$.
For $0<\kappa\leq\kappa_*$, the domain of attraction of $\nu_\kappa$ is
\[
    \mathcal D_\kappa:=\left\{
        \mu\in\mathcal P((0,\infty)):
        \mbb P_\mu(Z_t\in\cdot\mid Z_t>0)
        \Longrightarrow\nu_\kappa
        \ \text{as }t\to\infty
      \right\}.
\]
We use the convention $-\log0=+\infty$ in the tail criterion below.

\bgtheorem[Complete domains of attraction]
\label{thm_complete_domains_attraction}
The domains of attraction are given by
\begin{itemize}
    \item[\rm(i)] For each $0<\kappa<\kappa_*$,
    \[
        \mathcal D_\kappa
        =\left\{
            \mu\in\mathcal P((0,\infty)):
            \lim_{x\to\infty}
            \frac{-\log\mu([x,\infty))}{\log x}=\kappa
          \right\}.
    \]
    \item[\rm(ii)] For the Yaglom distribution,
    \[
        \mathcal D_{\kappa_*}
        =\left\{
            \mu\in\mathcal P((0,\infty)):
            \int_0^\infty x^p\,\mu(\d x)<\infty
            \ \text{for every }0<p<\kappa_*
          \right\}.
    \]
\end{itemize}
\edtheorem

These criteria determine, for every initial law, whether it belongs to a
QSD domain of attraction and, if so, to which one.
{The conditional law from an initial distribution outside these domains
does not converge to any QSD.}
{In particular, every initial law with finite first moment belongs to
$\mathcal D_{\kappa_*}$, since $\kappa_*\leq1$.}

\begin{rem1}[Tails without regular variation]\rm{}
\label{rem_oscillating_tail}
{Several attraction criteria for branching processes and diffusions
rely on regular variation assumptions; see \cite{RVJ68,LSM00,Yam22}.
For the CBBRE, Theorem \ref{thm_complete_domains_attraction}(i) characterizes
attraction through a logarithmic tail condition that is strictly weaker
than regular variation of the tail.}

{To illustrate the difference, fix} $0<\kappa<\kappa_*$ and
$0<\varepsilon<\kappa/\sqrt{\kappa^2+1}$.
The formula
\[
    {\mu_0}([x,\infty))
    =x^{-\kappa}\bigl(1+\varepsilon\sin(\log x)\bigr),
    \qquad x\geq1,
\]
defines a probability law supported on $[1,\infty)$. Indeed, the right-hand
side equals $1$ at $x=1$, tends to zero, and has derivative at most
$x^{-\kappa-1}(-\kappa+\varepsilon\sqrt{\kappa^2+1})<0$.
Its logarithmic tail exponent is $\kappa$, so
Theorem \ref{thm_complete_domains_attraction}(i) gives
${\mu_0}\in\mathcal D_\kappa$. However,
\[
    \frac{{\mu_0}([\e^\pi x,\infty))}{{\mu_0}([x,\infty))}
    =\e^{-\kappa\pi}
      \frac{1-\varepsilon\sin(\log x)}
           {1+\varepsilon\sin(\log x)}
\]
has no limit as $x\to\infty$. Thus regular variation of the tail is not
necessary for attraction to an interior QSD.
\end{rem1}

\subsection{Proof strategy}\label{sec_intro_proof_strategy}

The two parts of the problem require different arguments.  For the
classification of QSDs, the eigenmeasure relation leads to a second-order
differential equation for the modified Laplace transform.  The equation can be
solved explicitly in terms of Kummer functions.
{However, its solutions need not be modified Laplace transforms of
probability measures. By examining the equation near the origin, we show
that this requirement forces $0<\kappa\leq\kappa_*$. To establish existence
throughout this range, we construct the interior QSDs $\nu_\kappa$ with
$0<\kappa<\kappa_*$
as conditional limits from initial laws with exact power tails.
The endpoint QSD $\nu_{\kappa_*}$ is the known Yaglom distribution.}

For the attraction problem, a key step is to determine which initial states
account for survival at large times.  For sufficiency, suppose first that the
logarithmic tail exponent of $\mu$ is $\kappa\in(0,\kappa_*)$.
{Our estimates show that the conditional law from $\mu$ is
asymptotically determined by initial states between two exponential levels.}
More precisely, with $v=-{\bf m}-\kappa\sigma^2$, for every
sufficiently small $\varepsilon>0$,
\begin{equation}\label{intro_initial_concentration}
 \lim_{t\to\infty}
 \mbb P_\mu\bigl(\e^{(v-\varepsilon)t}\le Z_0
              \le\e^{(v+\varepsilon)t}\mid Z_t>0\bigr)=1.
\end{equation}
{This follows from Lemma \ref{full_lem_tail_survival_equivalence} and
Proposition \ref{full_prop_interior_concentration}.}

To see how \eqref{intro_initial_concentration} determines the conditional
limit, write, for each fixed $\lambda>0$,
\begin{equation}\label{intro_conditional_mixture}
 \mbb E_\mu[1-\e^{-\lambda Z_t}\mid Z_t>0]
 =\int_{(0,\infty)}
   \mbb E_z[1-\e^{-\lambda Z_t}\mid Z_t>0]
   \,\mbb P_\mu(Z_0\in\d z\mid Z_t>0).
\end{equation}
By \eqref{intro_initial_concentration}, we may restrict the integral in
\eqref{intro_conditional_mixture} to initial states between
$\e^{(v-\varepsilon)t}$ and $\e^{(v+\varepsilon)t}$, with an error tending
to zero.  Lemma \ref{full_lem_moving_initial_mass} gives, for $\lambda>0$,
\begin{equation}\label{intro_moving_initial_limit}
 \begin{aligned}
 \mbb E_{\e^{vt}}[1-\e^{-\lambda Z_t}\mid Z_t>0]
 &\longrightarrow \Phi_{(-{\bf m}-v)/\sigma^2}(\lambda)
 =\Phi_\kappa(\lambda)\\
 &=\int_{(0,\infty)}(1-\e^{-\lambda z})\,\nu_\kappa(\d z),
 \qquad\text{as }t\to\infty.
 \end{aligned}
\end{equation}
Using \eqref{intro_initial_concentration} and monotonicity (Lemma
\ref{full_lem_initial_mass_monotonicity}), we apply
\eqref{intro_moving_initial_limit} in \eqref{intro_conditional_mixture}
to obtain
\[
 \lim_{t\to\infty}
 \mbb E_\mu[1-\e^{-\lambda Z_t}\mid Z_t>0]
 =\Phi_\kappa(\lambda),\qquad\lambda>0.
\]
{Thus the conditional law from $\mu$ converges weakly to $\nu_\kappa$.}

For necessity, suppose that {the conditional law from $\mu$ converges to $\nu_\kappa$}
with $0<\kappa<\kappa_*$.  This convergence fixes the exponential survival
rate at $\gamma(\kappa)$.
{The key point is that this survival rate forces the logarithmic tail
exponent of $\mu$ to converge to $\kappa$.
We first deduce that every positive moment of order below $\kappa$ is
finite. Markov's inequality then gives
\begin{equation}\label{intro_survival_tail_rate}
 \liminf_{x\to\infty}\frac{-\log\mu([x,\infty))}{\log x}\geq\kappa.
\end{equation}
See Lemmas \ref{full_lem_conditional_survival_rate} and
\ref{full_lem_survival_rates}(ii).}

{
A heavier initial tail makes large starting populations more likely and can
therefore slow the exponential decay of survival.  Thus
\eqref{intro_survival_tail_rate} rules out logarithmic tail decay at a
rate below $\kappa$, but still allows faster decay along some sequences.
The next estimate rules out that possibility.}

{The survival rate allows us to apply Proposition
\ref{full_prop_interior_concentration}, so
\eqref{intro_initial_concentration} holds here as well.} We use this estimate
to bound the initial tail by a ratio of survival probabilities.
Combining the known survival
rate $\gamma(\kappa)$ for the numerator with the survival asymptotics in
Lemma \ref{full_lem_moving_initial_mass} for the denominator, we obtain
\begin{equation}\label{intro_necessary_tail_bound}
 \begin{aligned}
 \mu\bigl([\e^{(v-\varepsilon)t},\infty)\bigr)
 &\ge (1-o(1))
 \frac{\mbb P_\mu(Z_t>0)}
      {\mbb P_{\e^{(v+\varepsilon)t}}(Z_t>0)}\\
 &=\exp\!\left\{-\bigl(\kappa v+O(\varepsilon)\bigr)t+o(t)\right\}.
 \end{aligned}
\end{equation}
Taking logarithms with $x=\e^{(v-\varepsilon)t}$ gives
\[
 \limsup_{x\to\infty}\frac{-\log\mu([x,\infty))}{\log x}
 \le \frac{\kappa v+O(\varepsilon)}{v-\varepsilon}
 \longrightarrow\kappa
 \qquad\text{as }\varepsilon\downarrow0.
\]
Together with \eqref{intro_survival_tail_rate}, this proves the required
tail limit.

The case $\kappa=\kappa_*$ is treated separately.  We use the moment
condition to bound the contribution from large initial states and combine
fixed-state Yaglom convergence with the limits for exponentially growing
initial states to prove sufficiency.
{For necessity, conditional convergence determines the survival
rate, which implies the moment condition by Lemma
\ref{full_lem_survival_rates}(ii).}

We now describe the main tools behind these arguments.
We express survival probabilities and conditional modified Laplace
transforms in terms of Brownian exponential functionals.
{Conditioning on the Brownian endpoint and using a Brownian-bridge
representation give the conditional limit \eqref{intro_moving_initial_limit}
and the survival asymptotics used in \eqref{intro_necessary_tail_bound}.
An auxiliary diffusion and a coupling argument establish monotonicity of the
conditional modified Laplace transform in the initial state; together with
\eqref{intro_initial_concentration}, this allows us to extend
\eqref{intro_moving_initial_limit} to general initial laws.
Moment bounds and survival estimates yield the concentration estimate
\eqref{intro_initial_concentration} and, under conditional convergence,
give the necessary tail bound \eqref{intro_survival_tail_rate}.}

The remainder of the paper is organized as follows.
Section \ref{sec_qsd_equation} derives the equation for the modified Laplace
transforms of the QSDs and proves Theorem \ref{thm_complete_qsd_classification}(i).
Section \ref{sec_exact_power_domains} establishes the moment estimates
needed to construct the QSDs and proves Theorem
\ref{thm_complete_qsd_classification}(ii).
Section \ref{sec_complete_domains_attraction} proves the complete attraction
criterion in Theorem \ref{thm_complete_domains_attraction}, using the
auxiliary SDE and moving initial masses.

\section{The QSD Equation for Laplace Transforms}\label{sec_qsd_equation}
\setcounter{equation}{0}

This section proves Theorem \ref{thm_complete_qsd_classification}(i) by
identifying the possible QSD transforms and decay rates.

For a bounded measurable function $g$ on $[0,\infty)$, write
$P_tg(z)=\mathbb E_z[g(Z_t)]$, $z,t\geq0$, for the transition semigroup.
Let $\nu$ be a QSD with decay rate $\rho>0$.  By
\eqref{qsd_left_eigenmeasure},
$\nu P_t^0=\e^{-\rho t}\nu$.  Since $0$ is absorbing,
$P_t^0g=P_tg$ on $(0,\infty)$ whenever $g(0)=0$.  Therefore, for every bounded
measurable function $g$ with $g(0)=0$,
\beqlb\label{qsd_semigroup}
    \int_0^\infty P_tg(z)\nu(\d z)
    =
    \e^{-\rho t}\int_0^\infty g(z)\nu(\d z).
\eeqlb
For $\lambda>0$, set $f_\lambda(z)=1-\e^{-\lambda z}$.
By \eqref{def_modified_laplace_general},
\beqlb\label{def_modified_laplace_qsd}
    H_\nu(\lambda)
    =
    \int_0^\infty f_\lambda(z)\nu(\d z).
\eeqlb
For $\lambda>0$, set
\beqlb\label{Lf_explicit}
b_\lambda(z)
\ar=\ar
-z\psi(\lambda)\e^{-\lambda z}
-\frac12\sigma^2\lambda^2z^2\e^{-\lambda z}\cr
\ar=\ar
\left(
    \alpha\lambda-c\lambda^{1+\beta}
    -\frac12\sigma^2\lambda^2z
\right)z\e^{-\lambda z}.
\eeqlb

\bgproposition\label{prop_weak_eigen}
For every $\lambda>0$, define
\[
    M_t^{(\lambda)}
    :=f_\lambda(Z_t)-f_\lambda(Z_0)
    -\int_0^t b_\lambda(Z_s)\,\d s,
    \qquad t\geq0,
\]
Then $(M_t^{(\lambda)})_{t\geq0}$ is a martingale under $\mbb P_z$ for every
$z\geq0$.  Moreover,
\beqlb\label{weak_eigen}
    \int_0^\infty b_\lambda(z)\nu(\d z)
    =
    -\rho H_\nu(\lambda).
\eeqlb
\edproposition

\proof
Applying It\^o's formula to \eqref{main_sde_feller} or \eqref{main_sde} gives
\beqlb\label{ito_laplace}
    f_\lambda(Z_t)
    =
    f_\lambda(Z_0)
    +\int_0^t b_\lambda(Z_s)\,\d s
    +M_t^{(\lambda)},
\eeqlb
where the stochastic integral representation of $M_t^{(\lambda)}$ depends on
$\beta$.  When $\beta=1$,
\[
    M_t^{(\lambda)}
    =
    \sigma\lambda\int_0^t
    Z_s\e^{-\lambda Z_s}\,\d B_s^{(e)}
    +
    \sqrt{2c}\lambda\int_0^t
    \sqrt{Z_s}\e^{-\lambda Z_s}\,\d B_s^{(b)}.
\]
The process $(M_t^{(\lambda)})_{t\geq0}$ is a continuous local martingale.

Suppose next that $0<\beta<1$.  The identity
\[
    \frac{c\beta(\beta+1)}{\Gamma(1-\beta)}
    \int_0^\infty
    \bigl(\e^{-\lambda y}-1+\lambda y\bigr)y^{-2-\beta}\,\d y
    =
    c\lambda^{1+\beta}
\]
and the jump version of It\^o's formula give
\[
 \begin{aligned}
    M_t^{(\lambda)}
    ={}&
    \sigma\lambda\int_0^t
    Z_s\e^{-\lambda Z_s}\,\d B_s^{(e)}\\
    &+
    \int_0^t\int_0^\infty\int_0^{Z_{s-}}
    \e^{-\lambda Z_{s-}}
    \bigl(1-\e^{-\lambda y}\bigr)
    \,\widetilde N(\d s,\d y,\d u).
 \end{aligned}
\]
In this case, $(M_t^{(\lambda)})_{t\geq0}$ is also a local martingale.

In both cases, $f_\lambda$ and $b_\lambda$ are bounded, with
$b_\lambda(0)=0$.  For every $T>0$, \eqref{ito_laplace} gives
\[
    \sup_{0\leq t\leq T}|M_t^{(\lambda)}|
    \leq 2\|f_\lambda\|_\infty+T\|b_\lambda\|_\infty.
\]
Applying dominated convergence along a localizing sequence shows that
$(M_t^{(\lambda)})_{t\geq0}$ is a martingale.
Taking expectations in \eqref{ito_laplace} and using boundedness of
$b_\lambda$ to apply Fubini's theorem gives
\beqlb\label{dynkin}
    P_tf_\lambda(z)-f_\lambda(z)
    =
    \int_0^t P_sb_\lambda(z)\,\d s.
\eeqlb
Integrating \eqref{dynkin} against $\nu$ and using
\eqref{qsd_semigroup}, we obtain
\[
    \bigl(\e^{-\rho t}-1\bigr)H_\nu(\lambda)
    =
    \left(\int_0^t\e^{-\rho s}\,\d s\right)
    \int_0^\infty b_\lambda(z)\nu(\d z).
\]
Applying the identity
\[
    \e^{-\rho t}-1
    =
    -\rho\int_0^t\e^{-\rho s}\,\d s
\]
gives \eqref{weak_eigen}. 
\qed

For $\rho\geq0$, consider the differential equation
\beqlb\label{master_ode}
\frac{1}{2}\sigma^2\lambda^2 H''(\lambda)
+(\alpha\lambda-c\lambda^{1+\beta})H'(\lambda)
+\rho H(\lambda)=0,\qquad\lambda>0.
\eeqlb

\bglemma\label{lem_generator_laplace}
Let $\nu$ be a quasi-stationary distribution satisfying
\eqref{qsd_semigroup} with decay rate $\rho$.  Its modified Laplace transform
$H_\nu$, given by \eqref{def_modified_laplace_qsd}, satisfies
\eqref{master_ode}.
Furthermore, $0\leq H_\nu\leq1$,
$\lim_{\lambda\downarrow0}H_\nu(\lambda)=0$, and
$\lim_{\lambda\to\infty}H_\nu(\lambda)=1$.
\edlemma

\proof
 Since $\nu$ is a probability measure, dominated
convergence justifies differentiating \eqref{def_modified_laplace_qsd}
twice under the integral.  Thus $H_\nu\in C^2((0,\infty))$ and, for every
$\lambda>0$,
\[
    H_\nu'(\lambda)
    =
    \int_0^\infty z\e^{-\lambda z}\nu(\d z),
    \qquad
    H_\nu''(\lambda)
    =
    -\int_0^\infty z^2\e^{-\lambda z}\nu(\d z).
\]
Substituting \eqref{Lf_explicit} into \eqref{weak_eigen} gives
\[
 \begin{aligned}
    -\rho H_\nu(\lambda)
    &=\int_0^\infty b_\lambda(z)\nu(\d z)\\
    &=(\alpha\lambda-c\lambda^{1+\beta})
      \int_0^\infty z\e^{-\lambda z}\nu(\d z)\\
    &\quad-\frac12\sigma^2\lambda^2
      \int_0^\infty z^2\e^{-\lambda z}\nu(\d z)\\
    &=(\alpha\lambda-c\lambda^{1+\beta})H_\nu'(\lambda)
      +\frac12\sigma^2\lambda^2H_\nu''(\lambda).
 \end{aligned}
\]
The bounds and boundary limits follow directly from
\eqref{def_modified_laplace_qsd} and dominated convergence.
\qed

We next solve the boundary-value problem associated with
\eqref{master_ode}.  For $0<\rho\leq\gamma(\kappa_c)$, consider
\beqlb\label{characteristic_eq}
    \frac12\beta^2\sigma^2p^2+\beta{\bf m}p+\rho=0.
\eeqlb
Its discriminant is $\beta^2({\bf m}^2-2\sigma^2\rho)\geq0$.
Denote its two real roots by $p_0\leq p_1$.  By Vieta's formulas,
\beqlb\label{kummer_root_relations}
    p_0+p_1=-\frac{2{\bf m}}{\beta\sigma^2}>0,
    \qquad
    p_0p_1=\frac{2\rho}{\beta^2\sigma^2}>0.
\eeqlb
Thus $0<p_0\leq p_1$, with equality if and only if
$\rho=\gamma(\kappa_c)$.  Set
\beqlb\label{kummer_parameters}
    c_0=\frac{2c}{\beta\sigma^2},
    \qquad
    b_0=2p_0+1+\frac{2{\bf m}}{\beta\sigma^2}=1+p_0-p_1.
\eeqlb
The last equality follows from \eqref{kummer_root_relations}.

\bglemma
\label{lem_ode_uniqueness}
For $0<\rho\leq\gamma(\kappa_c)$, equation \eqref{master_ode} has a unique
solution $H\in C^2((0,\infty))$ with $H(0+)=0$ and $H(\infty)=1$, namely
\beqlb\label{pure_ode_solution}
H(\lambda)
=
(c_0\lambda^\beta)^{p_0}U(p_0,b_0,c_0\lambda^\beta).
\eeqlb
\edlemma

\proof
With the constants in \eqref{kummer_parameters}, put
$u=c_0\lambda^\beta$ and write $H(\lambda)=u^{p_0}W(u)$.
The chain rule and \eqref{characteristic_eq} transform
\eqref{master_ode} into Kummer's equation
\beqlb\label{kummer_equation_for_W}
    uW''(u)+(b_0-u)W'(u)-p_0W(u)=0.
\eeqlb

The function $U(p_0,b_0,\cdot)$ solves \eqref{kummer_equation_for_W} and
satisfies $U(p_0,b_0,u)\sim u^{-p_0}$ as $u\to\infty$; see
\cite[(13.2.1) and (13.2.6)]{NISTDLMF}.
The fundamental pair in \cite[(13.2.24)]{NISTDLMF}, together with the
asymptotic formula \cite[(13.2.6)]{NISTDLMF}, shows that every solution
$W$ of \eqref{kummer_equation_for_W} linearly independent of
$U(p_0,b_0,\cdot)$ satisfies $W(u)\sim C\e^u u^{p_0-b_0}$ as
$u\to\infty$, for some $C\neq0$.
Thus the condition $\lim_{\lambda\to\infty}H(\lambda)=1$ uniquely
determines $W(u)=U(p_0,b_0,u)$.

It remains to verify that the candidate satisfies $H(0+)=0$.
The integral representation \cite[(13.4.4)]{NISTDLMF}, with the
parameters in \eqref{kummer_parameters}, followed by the substitution
$v=ur$, gives
\[
\begin{aligned}
    H(\lambda)
    &=\frac{u^{p_0}}{\Gamma(p_0)}
      \int_0^\infty \e^{-ur}r^{p_0-1}(1+r)^{-p_1}\,\d r\\
    &=\frac{1}{\Gamma(p_0)}
      \int_0^\infty \e^{-v}v^{p_0-1}
      (1+v/u)^{-p_1}\,\d v,
      \qquad u=c_0\lambda^\beta.
\end{aligned}
\]
For every $v>0$, the factor $(1+v/u)^{-p_1}$ tends to zero as
$u\downarrow0$ and is bounded by $1$.
Since $\e^{-v}v^{p_0-1}$ is integrable, dominated convergence yields
$H(0+)=0$.  This completes the proof of existence and uniqueness for
the boundary-value problem.

Finally, the expression does not depend on which root is used.
Kummer's transformation \cite[(13.2.40)]{NISTDLMF}, together with
\eqref{kummer_parameters}, gives
\[
\begin{aligned}
    u^{p_0}U(p_0,b_0,u)
    &=u^{p_0+1-b_0}U(p_0-b_0+1,2-b_0,u)\\
    &=u^{p_1}U(p_1,1+p_1-p_0,u).
\end{aligned}
\]
Thus changing the root changes both the prefactor and the parameters
of $U$, but leaves the resulting function $H$ unchanged.
\qed

\begin{cor1}\label{main_thm_1}
Let $\nu$ be a QSD for the CBBRE with decay rate
$0<\rho\leq\gamma(\kappa_c)$.
Then its modified Laplace transform is given by \eqref{pure_ode_solution}.
Consequently, there is at most one QSD with each such decay rate.
\end{cor1}

\proof
By Lemma \ref{lem_generator_laplace}, $H_\nu$ satisfies
\eqref{master_ode} with $H_\nu(0+)=0$ and $H_\nu(\infty)=1$.
Lemma \ref{lem_ode_uniqueness} gives $H_\nu$ in the form
\eqref{pure_ode_solution}, and
uniqueness of Laplace transforms gives uniqueness of the QSD.
\qed

Lemma \ref{lem_ode_uniqueness} solves the boundary-value problem for
every $0<\rho\leq\gamma(\kappa_c)$.  Such a solution need not be the
modified Laplace transform of a probability measure.  The next lemma
shows that every QSD has decay rate at most $\gamma(\kappa_*)$.
In particular, when $\kappa_*<\kappa_c$,
no QSD has decay rate in $(\gamma(\kappa_*),\gamma(\kappa_c)]$.

\bglemma\label{lem_origin_admissibility}
Let $\nu$ be a QSD for the CBBRE with decay rate $\rho>0$.
Then $\rho\leq\gamma(\kappa_*)$.
\edlemma

\proof
Suppose on the contrary that $\rho>\gamma(\kappa_*)$ and put
$\delta=\rho-\gamma(\kappa_*)>0$.
For $\lambda>0$, define
$R(\lambda):=\lambda H_\nu'(\lambda)/H_\nu(\lambda)$.
The integral representation \eqref{def_modified_laplace_qsd} and the
inequality $0<u\e^{-u}\leq1-\e^{-u}$ for $u>0$ give
\beqlb\label{qsd_ratio_bound}
    0<R(\lambda)\leq1,\qquad \lambda>0.
\eeqlb

By Lemma \ref{lem_generator_laplace}, $H_\nu$ satisfies
\eqref{master_ode}.  Differentiating $R$ and using this equation gives
\beqlb\label{qsd_ratio_ode}
\begin{aligned}
    \lambda R'(\lambda)
    &=R(\lambda)+\frac{\lambda^2H_\nu''(\lambda)}{H_\nu(\lambda)}
      -R(\lambda)^2\\
    &=\frac{2}{\sigma^2}
      \left\{\gamma(R(\lambda))-\rho+c\lambda^\beta R(\lambda)\right\}.
\end{aligned}
\eeqlb
Here we used $\alpha-\sigma^2/2={\bf m}$ and
\eqref{gamma_kappa_front}.
Since $\gamma'(r)=\sigma^2(\kappa_c-r)$ and
$\kappa_*=\min\{1,\kappa_c\}$, we have
$\gamma(r)\leq\gamma(\kappa_*)$ for $0\leq r\leq1$.

Choose $\lambda_0>0$ small enough that $c\lambda_0^\beta\leq\delta/2$.
Combining \eqref{qsd_ratio_bound} and \eqref{qsd_ratio_ode} gives
\[
    \lambda R'(\lambda)
    \leq\frac{2}{\sigma^2}(-\delta+c\lambda^\beta)
    \leq-\frac{\delta}{\sigma^2},
    \qquad 0<\lambda\leq\lambda_0.
\]
Dividing by $\lambda$ and integrating from $\lambda$ to $\lambda_0$,
we obtain
\[
    R(\lambda)
    \geq R(\lambda_0)+\frac{\delta}{\sigma^2}
        \log\frac{\lambda_0}{\lambda}
    \longrightarrow\infty
    \qquad\text{as }\lambda\downarrow0.
\]
This contradicts \eqref{qsd_ratio_bound}.
\qed

\noindent{\it Proof of Theorem \ref{thm_complete_qsd_classification}(i).~}
Let $\nu$ be a QSD with decay rate $\rho$.
Lemma \ref{lem_origin_admissibility} gives
$0<\rho\leq\gamma(\kappa_*)$.
Since $\gamma(0)=0$ and $\gamma$ is continuous and strictly increasing
on $[0,\kappa_*]$, there is a unique $\kappa$ such that
\[
    \rho=\gamma(\kappa),
    \qquad 0<\kappa\leq\kappa_*.
\]
Corollary \ref{main_thm_1}, with $p_0=\kappa/\beta$, then gives
$H_\nu=\Phi_\kappa$ and shows that there is at most one QSD with each
decay rate.
\qed

The existence assertion in Theorem \ref{thm_complete_qsd_classification}(ii)
will be proved in the following section.

\section{Construction of quasi-stationary distributions}\label{sec_exact_power_domains}
\setcounter{equation}{0}

We prove Theorem \ref{thm_complete_qsd_classification}(ii).
We first derive moment estimates for the conditional Laplace exponent,
then use them to construct the interior QSDs. The Yaglom distribution
provides the remaining member.  The moment estimates will also be used
in Section \ref{sec_complete_domains_attraction}.

We begin with a criterion for a conditional limit to be a QSD.

\bglemma\label{lem_conditional_limit_qsd}
Let $\mu$ be an initial probability law on $(0,\infty)$.  Suppose that, for
some probability law $\nu$ on $(0,\infty)$,
\beqlb\label{conditional_limit_weak_convergence}
    \mbb P_\mu(Z_t\in\cdot\mid Z_t>0)
    \Longrightarrow \nu
    \qquad\text{as }t\to\infty.
\eeqlb
Then the following assertions hold.
\begin{itemize}
    \item[\rm(i)] The probability measure $\nu$ is a QSD.
    \item[\rm(ii)] If, for some $\rho\geq0$ and every fixed $s\geq0$,
    \beqlb\label{conditional_limit_survival_ratio}
        \frac{\mbb P_\mu(Z_{t+s}>0)}
             {\mbb P_\mu(Z_t>0)}
        \xrightarrow[t\to\infty]{}
        \e^{-\rho s},
    \eeqlb
    then $\nu P_s^0=\e^{-\rho s}\nu$ for every $s\geq0$.  In particular,
    the decay rate of $\nu$ is $\rho$.
\end{itemize}
\edlemma

\proof
Write $\mu_t=\mbb P_\mu(Z_t\in\cdot\mid Z_t>0)$ and
$e_q(z)=\e^{-qz}$.  Put $S_r=\sigma B_r^{(e)}+{\bf m}r$ and, for
$s,q>0$, define
\[
    V_s(q)=\left(q^{-\beta}\e^{-\beta S_s}
        +c\beta\int_0^s\e^{-\beta S_r}\,\d r\right)^{-1/\beta}.
\]
Set $V_s(0)=0$ and $V_s(\infty)=\lim_{q\to\infty}V_s(q)$.
By \cite[Example~2, (14), and Example~3, (16)]{PP17}, the conditional
Laplace transform given the environment is
\[
    \mbb E_z[\e^{-qZ_s}\mid (S_r)_{0\leq r\leq s}]
    =\e^{-zV_s(q)},\qquad z,q\geq0.
\]
In particular, $V_s(q)$ does not depend on $z$, and
$0<V_s(\infty)<\infty$ almost surely.  Taking expectations and letting
$q\to\infty$, bounded convergence gives
$\mbb P_z(Z_s=0)=\mbb E[\e^{-zV_s(\infty)}]$.
Subtracting this extinction probability from the unconditional Laplace
transform yields, for $s>0$ and $q\geq0$,
\[
    P_s^0e_q(z)
    =
    \mbb E\!\left[
        \e^{-zV_s(q)}-\e^{-zV_s(\infty)}
    \right],
    \qquad
    P_s^0\mathbf1(z)
    =
    \mbb E[1-\e^{-zV_s(\infty)}].
\]
Both functions are bounded and continuous in $z$, and
$P_s^0\mathbf1(z)>0$ for $z>0$.  By the Markov property and absorption
at zero, for $g=e_q$ or $g=\mathbf1$,
\beqlb\label{conditional_limit_markov_expectation}
\begin{aligned}
    \mbb E_\mu\!\left[g(Z_{t+s})\mathbf1_{\{Z_{t+s}>0\}}\right]
    &=\mbb E_\mu\!\left[
        \mathbf1_{\{Z_t>0\}}P_s^0g(Z_t)\right]\\
    &=\mbb P_\mu(Z_t>0)\,\mu_t(P_s^0g).
\end{aligned}
\eeqlb
Consequently,
\[
    \mu_{t+s}(e_q)
    =
    \frac{\mbb E_\mu\!\left[
        e_q(Z_{t+s})\mathbf1_{\{Z_{t+s}>0\}}\right]}
         {\mbb P_\mu(Z_{t+s}>0)}
    =
    \frac{\mu_t(P_s^0e_q)}{\mu_t(P_s^0\mathbf1)}.
\]
By \eqref{conditional_limit_weak_convergence}, both $\mu_t$ and
$\mu_{t+s}$ converge weakly to $\nu$.  Letting $t\to\infty$ gives
\beqlb\label{conditional_limit_laplace_identity}
    \nu(e_q)
    =
    \frac{\nu(P_s^0e_q)}{\nu(P_s^0\mathbf1)}.
\eeqlb
Note that $\nu P_s^0/\nu(P_s^0\mathbf1)$ is the conditional law
$\mbb P_\nu(Z_s\in\cdot\mid Z_s>0)$.
Equation \eqref{conditional_limit_laplace_identity} shows that this law
has the same Laplace transform as $\nu$, so they are equal by
uniqueness of Laplace transforms.  This proves part~(i).

For part~(ii), taking $g=\mathbf1$ in
\eqref{conditional_limit_markov_expectation} gives
\beqlb\label{conditional_limit_survival_identity}
    \mu_t(P_s^0\mathbf1)
    =
    \frac{\mbb P_\mu(Z_{t+s}>0)}
         {\mbb P_\mu(Z_t>0)}.
\eeqlb
By \eqref{conditional_limit_weak_convergence} and the bounded continuity
of $P_s^0\mathbf1$, we may let $t\to\infty$ in
\eqref{conditional_limit_survival_identity}.  Using
\eqref{conditional_limit_survival_ratio} yields
\beqlb\label{conditional_limit_survival_limit}
    \nu(P_s^0\mathbf1)=\e^{-\rho s}.
\eeqlb
Substituting \eqref{conditional_limit_survival_limit} into
\eqref{conditional_limit_laplace_identity} gives
$\nu(P_s^0e_q)=\e^{-\rho s}\nu(e_q)$ for every $q\geq0$.
Uniqueness of Laplace transforms yields $\nu P_s^0=\e^{-\rho s}\nu$.
\qed

\subsection{Moment asymptotics}\label{sec_strict_fractional_moments}

To construct QSDs through conditional limits, we estimate the conditional
Laplace transform.  For an initial probability law $\mu$
and $0<\theta<1$, Karamata's Tauberian theorem gives
\beqlb\label{exact_power_tauberian_equivalence}
    \mu((x,\infty))\sim c_\mu x^{-\theta}
    \quad (x\to\infty)
    \quad\Longleftrightarrow\quad
    H_\mu(r)\sim c_\mu\Gamma(1-\theta)r^\theta
    \quad (r\downarrow0);
\eeqlb
see \cite[Corollary~8.1.7]{BGT87}.  If $\mu$ has finite first moment, monotone
convergence instead gives
\[
    H_\mu(r)\sim r\int_0^\infty z\,\mu(\d z),
    \qquad r\downarrow0.
\]
These formulas translate the initial-tail assumptions into the small-$r$
asymptotics used below.
Recall that $S_t=\sigma B_t^{(e)}+{\bf m}t$ for $t\geq0$.
For $t>0$ and $\lambda\in(0,\infty]$, write
\beqlb\label{Vt_environment_definition}
    V_t(\lambda)
    :=\left(
        \lambda^{-\beta}\e^{-\beta S_t}
        +c\beta\int_0^t\e^{-\beta S_r}\,\d r
    \right)^{-1/\beta},
\eeqlb
where $\infty^{-\beta}=0$.
Using the conditional Laplace formula of Palau and Pardo
\cite[(14) and (16)]{PP17}, we obtain
\beqlb\label{annealed_ratio}
    \mbb E_\mu[1-\e^{-\lambda Z_t}\mid Z_t>0]
    =
    \frac{\mbb E[H_\mu(V_t(\lambda))]}
         {\mbb E[H_\mu(V_t(\infty))]},
    \qquad t,\lambda>0.
\eeqlb
The next lemma expresses $V_t(\lambda)$ in terms of a Brownian
exponential functional.  This is the representation used in
\cite[(2.1) and Theorem~2.1]{LZZ}; we give the scaling argument to
make the notation explicit.

\bglemma
\label{lem_brownian_exponential_representation}
Set
\beqlb\label{brownian_scaling_parameters}
    \tau_t:=\frac{\sigma^2\beta^2}{4}t,\qquad
    \eta:=-\frac{2{\bf m}}{\beta\sigma^2},\qquad
    \gamma_0:=\frac{4c}{\beta\sigma^2}.
\eeqlb
Define
\beqlb\label{rescaled_environment_brownian}
    B_u:=-\frac{\sigma\beta}{2}
        B^{(e)}_{4u/(\sigma^2\beta^2)},\qquad u\geq0,
\eeqlb
and, for $h\in\mbb R$, put
\[
    A_s^{(h)}:=\int_0^s \e^{2(B_u+hu)}\,\d u,\qquad s\geq0.
\]
Then $(B_u)_{u\geq0}$ is a standard Brownian motion and, for
$t>0$ and $\lambda\in(0,\infty]$,
\beqlb\label{Vt_explicit_representation}
    V_t(\lambda)
    =
    \left(
        \lambda^{-\beta}\e^{2(B_{\tau_t}+\eta\tau_t)}
        +\gamma_0 A_{\tau_t}^{(\eta)}
    \right)^{-1/\beta},
\eeqlb
where $\infty^{-\beta}=0$.  In particular,
$V_t(\infty)=(\gamma_0 A_{\tau_t}^{(\eta)})^{-1/\beta}$.
Equivalently, for a standard Brownian motion $B$,
\begin{equation}\label{full_eq_physical_representation}
    V_t(\lambda)\stackrel{\d}{=}
    \e^{\sigma B_t+{\bf m}t}
    \left(\lambda^{-\beta}+c\beta\int_0^t
          \e^{\beta(\sigma B_s+{\bf m}s)}\,\d s\right)^{-1/\beta},
    \qquad \lambda\in(0,\infty].
\end{equation}
This is a fixed-time identity jointly in $\lambda$.
\edlemma

\proof
By Brownian scaling and symmetry
\cite[Lemma~2.9.4(i), (iv)]{KS91}, $(B_u)_{u\geq0}$ defined by
\eqref{rescaled_environment_brownian} is a standard Brownian motion
with respect to its natural filtration.
By \eqref{brownian_scaling_parameters},
\[
    -\beta(\sigma B_r^{(e)}+{\bf m}r)
    =2(B_{\tau_r}+\eta\tau_r),\qquad r\geq0.
\]
Substituting $u=\sigma^2\beta^2r/4$ therefore gives
\[
\begin{aligned}
    c\beta\int_0^t
        \e^{-\beta(\sigma B_r^{(e)}+{\bf m}r)}\,\d r
    &=\frac{4c}{\beta\sigma^2}
        \int_0^{\tau_t}\e^{2(B_u+\eta u)}\,\d u\\
    &=\gamma_0 A_{\tau_t}^{(\eta)}.
\end{aligned}
\]
Substituting these identities into \eqref{Vt_environment_definition} gives
\eqref{Vt_explicit_representation}.

For the equivalent representation, put $s=\tau_t$ and factor
$\e^{2(B_s+\eta s)}$ from the expression inside the parentheses in
\eqref{Vt_explicit_representation}. Time reversal and symmetry of Brownian
motion \cite[Lemma~2.9.4(iii), (iv)]{KS91} give
\[
    V_t(\lambda)\stackrel{\d}{=}
    \e^{2(B_s-\eta s)/\beta}
    \left(\lambda^{-\beta}+\gamma_0A_s^{(-\eta)}\right)^{-1/\beta}.
\]
By \eqref{brownian_scaling_parameters}, jointly in law,
\[
    \frac2\beta(B_{\tau_t}-\eta\tau_t)
    \stackrel{\d}{=}\sigma B_t+{\bf m}t,
    \qquad
    \gamma_0A_{\tau_t}^{(-\eta)}
    \stackrel{\d}{=}c\beta\int_0^t
          \e^{\beta(\sigma B_r+{\bf m}r)}\,\d r.
\]
This proves \eqref{full_eq_physical_representation}.
\qed

We use \eqref{Vt_explicit_representation} to estimate the moments of
$V_t(\lambda)$.

\bglemma\label{lem_fractional_moment_strict}
Let $p\in(0,\eta/2)$ and put $\kappa=\beta p$.  For
$\lambda\in(0,\infty]$, define
\[
    \mathcal M_{p,t}(\lambda)
    :=
    \mathbb E\!\left[
        \bigl(
            \lambda^{-\beta}\e^{2(B_{\tau_t}+\eta\tau_t)}
            +
            \gamma_0A_{\tau_t}^{(\eta)}
        \bigr)^{-p}
    \right],
\]
with the convention $\lambda^{-\beta}=0$ when $\lambda=\infty$.  There exists a
function $F_p:(0,\infty]\to(0,\infty)$ with
$F_p(\lambda)\leq\lambda^{\beta p}$ for $0<\lambda<\infty$ such that,
for every $\lambda\in(0,\infty]$,
\[
    \mathcal M_{p,t}(\lambda)
    \sim
    \e^{-\gamma(\kappa)t}F_p(\lambda),
    \qquad t\to\infty,
\]
where $\gamma$ is defined in \eqref{gamma_kappa_front}.
The limiting constant is
\begin{equation}\label{full_eq_fractional_constant}
    F_p(\lambda)
    =\mbb E\left[
       \left(\lambda^{-\beta}+c\beta\int_0^\infty
             \e^{\beta(\sigma B_s+({\bf m}+\kappa\sigma^2)s)}\,\d s\right)^{-p}
      \right],
    \,.
\end{equation}
Moreover, if $\kappa\leq\kappa_*$, then for every $\lambda\in(0,\infty)$,
\[
    \frac{F_p(\lambda)}{F_p(\infty)}
    =
    \Phi_\kappa(\lambda),
\]
where $\Phi_\kappa$ is defined in \eqref{def_Phi_kappa}.
\edlemma

\proof
Write $s=\tau_t$ and $a=\lambda^{-\beta}$. The time reversal in the
proof of Lemma \ref{lem_brownian_exponential_representation} gives
\[
    \mathcal M_{p,t}(\lambda)
    =
    \mathbb E\left[
        \e^{2p(B_s-\eta s)}
        \left(
            a+\gamma_0A_s^{(-\eta)}
        \right)^{-p}
    \right].
\]
On the fixed interval $[0,s]$, the constant integrand $2p$ satisfies
Novikov's condition, since
$\mbb E[\exp\{\frac12\int_0^s(2p)^2\,\d u\}]=\e^{2p^2s}<\infty$.
Thus $(\e^{2pB_u-2p^2u})_{0\leq u\leq s}$ is a martingale by
\cite[Corollary~3.5.13]{KS91}.  Under the probability measure with
density $\e^{2pB_s-2p^2s}$ relative to $\mbb P$ on
$\sigma(B_u:0\leq u\leq s)$, the process
$(B_u-2pu)_{0\leq u\leq s}$ is a standard Brownian motion by
\cite[Theorem~3.5.1]{KS91}.  Consequently,
\[
    \mathcal M_{p,t}(\lambda)
    =
    \e^{-2p(\eta-p)s}
    \mathbb E\left[
        \left(
            a+\gamma_0A_s^{(2p-\eta)}
        \right)^{-p}
    \right].
\]
Since $s=\sigma^2\beta^2t/4$ and $\kappa=\beta p$,
\[
    2p(\eta-p)s
    =
    \gamma(\kappa)t.
\]
The strict interior condition $p<\eta/2$ means that
\[
    \delta:=\eta-2p>0.
\]
Hence $A_s^{(2p-\eta)}=A_s^{(-\delta)}$ increases almost surely to
\[
    A_\infty^{(-\delta)}
    =
    \int_0^\infty \e^{2(B_u-\delta u)}\d u
    <\infty .
\]
Finiteness follows from the Brownian strong law
\cite[Problem~2.9.3]{KS91}: almost surely $B_u/u\to0$, so the
integrand is bounded by $\e^{-\delta u}$ for all sufficiently large $u$.
For finite $\lambda$, the integrand is bounded by $a^{-p}$, and dominated
convergence yields
\[
    F_p(\lambda)
    =
    \mathbb E\left[
        \left(
            \lambda^{-\beta}
            +
            \gamma_0A_\infty^{(-\delta)}
        \right)^{-p}
    \right].
\]
For $\lambda=\infty$, the same conclusion follows by dominated convergence:
for $s\geq s_0>0$, the integrand is dominated by a constant multiple of
$\bigl(A_{s_0}^{(-\delta)}\bigr)^{-p}$.  This negative moment is finite because
\[
    A_{s_0}^{(-\delta)}
    \geq
    s_0\exp\left\{2\inf_{0\leq u\leq s_0}(B_u-\delta u)\right\},
\]
and the reflection principle gives
$\sup_{0\leq u\leq s_0}(-B_u)\stackrel{d}{=}|B_{s_0}|$;
see \cite[Chapter~2, (8.3)]{KS91}.  Hence
\[
    \mbb E\!\left[\bigl(A_{s_0}^{(-\delta)}\bigr)^{-p}\right]
    \leq s_0^{-p}\e^{2p\delta s_0}
        \mbb E[\e^{2p|B_{s_0}|}]
    \leq 2s_0^{-p}\e^{2p\delta s_0+2p^2s_0}<\infty.
\]
The inverse time change in \eqref{brownian_scaling_parameters}, together
with $\delta=2(-{\bf m}-\kappa\sigma^2)/(\beta\sigma^2)$, turns the last representation of
$F_p$ into \eqref{full_eq_fractional_constant}.

For finite $\lambda$, the bound
$(\lambda^{-\beta}+\gamma_0A_\infty^{(-\delta)})^{-p}\leq\lambda^{\beta p}$
gives $F_p(\lambda)\leq\lambda^{\beta p}$.

It remains to identify the ratio.  Dufresne's identity \cite{Duf90}, in the
form given by \cite[Theorem~6.2]{MY05}, gives
\[
    A_\infty^{(-\delta)}
    \stackrel{d}{=}
    \frac1{2G_\delta},
\]
where $G_\delta$ has Gamma law with shape parameter $\delta$ and unit rate.  Hence
\[
    F_p(\infty)
    =
    \left(\frac2{\gamma_0}\right)^p
    \frac{\Gamma(\eta-p)}{\Gamma(\eta-2p)}.
\]
For finite $\lambda$, with
\[
    z=
    \frac{\gamma_0}{2}\lambda^\beta
    =
    \frac{2c}{\beta\sigma^2}\lambda^\beta,
\]
the integral representation \cite[(13.4.4)]{NISTDLMF} gives
\[
    \frac{F_p(\lambda)}{F_p(\infty)}
    =
    z^{\eta-p}U(\eta-p,1+\eta-2p,z).
\]
By Kummer's transformation \cite[(13.2.40)]{NISTDLMF}, this is equal to
\[
    z^pU(p,2p+1-\eta,z).
\]
Since
\[
    2p+1-\eta
    =
    2p+1+\frac{2{\bf m}}{\beta\sigma^2},
\]
for $\kappa\leq\kappa_*$, \eqref{def_Phi_kappa} identifies the last
expression with $\Phi_\kappa(\lambda)$.
\qed

\bglemma\label{lem_exact_power_transfer}
Assume
\[
    H_\mu(r)\sim C_\mu r^\kappa,
    \qquad r\downarrow0,
\]
with $C_\mu\in(0,\infty)$ and
\[
    0<\kappa<\kappa_c.
\]
Then, for every fixed $\lambda\in(0,\infty]$,
\[
    \mathbb E[H_\mu(V_t(\lambda))]
    \sim
    C_\mu\mathbb E[V_t(\lambda)^\kappa]
    \sim
    C_\mu \e^{-\gamma(\kappa)t}F_{\kappa/\beta}(\lambda),
    \qquad t\to\infty.
\]
\edlemma

\proof
Fix $\lambda\in(0,\infty]$ and define $(X_t)_{t>0}$ by
$X_t=V_t(\lambda)$.  For
$\varepsilon>0$, set
\[
    \omega(\varepsilon)
    =
    \sup_{0<r\leq\varepsilon}
    \left|
        \frac{H_\mu(r)}{C_\mu r^\kappa}-1
    \right|.
\]
Then $\omega(\varepsilon)\to0$ as $\varepsilon\downarrow0$, and
\[
    \left|
        \mathbb E[H_\mu(X_t);X_t\leq\varepsilon]
        -C_\mu\mathbb E[X_t^\kappa;X_t\leq\varepsilon]
    \right|
    \leq
    C_\mu\omega(\varepsilon)\mathbb E[X_t^\kappa].
\]

It remains to show that the complementary event is negligible on the
$\e^{-\gamma(\kappa)t}$ scale.  Choose
\[
    \kappa'\in(\kappa,\kappa_c).
\]
Since $\gamma$ is strictly increasing on $(0,\kappa_c)$,
\[
    \gamma(\kappa')>\gamma(\kappa).
\]
By Markov's inequality and Lemma \ref{lem_fractional_moment_strict},
\[
    \mathbb P(X_t>\varepsilon)
    \leq
    \varepsilon^{-\kappa'}\mathbb E[X_t^{\kappa'}]
    =
    O(\e^{-\gamma(\kappa')t})
    =
    o(\e^{-\gamma(\kappa)t}).
\]
Since $0\leq H_\mu\leq1$, the contribution of $\{X_t>\varepsilon\}$ is
negligible.  The same estimate also shows that
$\mathbb E[X_t^\kappa;X_t>\varepsilon]=o(\e^{-\gamma(\kappa)t})$.
Combining these estimates with Lemma \ref{lem_fractional_moment_strict}, first
letting $t\to\infty$ and then $\varepsilon\downarrow0$, proves the claim.
\qed

\subsection{Proof of Theorem \ref{thm_complete_qsd_classification}(ii)}
\label{subsec_qsd_construction}

We first construct the interior QSDs using the preceding moment estimates.

For an initial law $\mu$, write
$L_t^\mu(\lambda)=\mbb E_\mu[\e^{-\lambda Z_t}\mid Z_t>0]$.
For $0<\theta<1$, define the exact-power tail class
\[
    \mathcal T_\theta
    =
    \left\{
        \mu\in\mathcal P((0,\infty)):
        \mu((x,\infty))\sim c_\mu x^{-\theta},\quad
        c_\mu\in(0,\infty),\quad x\to\infty
    \right\}.
\]

\bglemma\label{lem_exact_power_construction}
For every $0<\kappa<\kappa_*$, there exists a QSD $\nu_\kappa$ with decay
rate $\gamma(\kappa)$ and modified Laplace transform $\Phi_\kappa$ such that,
for every $\mu\in\mathcal T_\kappa$,
\[
    \mbb P_\mu(Z_t\in\cdot\mid Z_t>0)
    \Longrightarrow
    \nu_\kappa.
\]
\edlemma

\proof
We first prove weak convergence.
If $\mu\in\mathcal T_\kappa$, then \eqref{exact_power_tauberian_equivalence} gives
\[
    H_\mu(r)\sim C_\mu r^\kappa,
    \qquad r\downarrow0,
\]
where $C_\mu=c_\mu\Gamma(1-\kappa)>0$.
By \eqref{annealed_ratio} and Lemmas \ref{lem_exact_power_transfer}
and \ref{lem_fractional_moment_strict}, for every $\lambda>0$,
\[
    \mbb E_\mu[1-\e^{-\lambda Z_t}\mid Z_t>0]
    =
    \frac{\mathbb E[H_\mu(V_t(\lambda))]}
         {\mathbb E[H_\mu(V_t(\infty))]}
    \to
    \frac{F_{\kappa/\beta}(\lambda)}
         {F_{\kappa/\beta}(\infty)}
    =
    \Phi_\kappa(\lambda),
    \qquad t\to\infty.
\]
Thus the conditional Laplace transforms converge to $1-\Phi_\kappa$.
By Lemma \ref{lem_fractional_moment_strict},
\[
    0\leq\Phi_\kappa(\lambda)
    \leq\frac{\lambda^\kappa}{F_{\kappa/\beta}(\infty)}
    \longrightarrow0,
    \qquad\lambda\downarrow0.
\]
Thus $1-\Phi_\kappa$ extends continuously to $0$ with value $1$.
By the continuity theorem for Laplace transforms
\cite[Chapter~XIII, Section~1, Theorem~2, p.~431]{Fel71}, the conditional
laws converge weakly to a probability measure $\nu_\kappa$ on $[0,\infty)$
with Laplace transform $1-\Phi_\kappa$.
Since $\Phi_\kappa(\lambda)\to1$ as $\lambda\to\infty$,
we have $\nu_\kappa(\{0\})=0$.

Lemma \ref{lem_conditional_limit_qsd}(i) now shows that $\nu_\kappa$ is a QSD.

It remains to identify its decay rate. Lemma \ref{lem_exact_power_transfer}
gives
\[
    \mathbb P_\mu(Z_t>0)
    =\mathbb E[H_\mu(V_t(\infty))]
    \sim C_\mu F_{\kappa/\beta}(\infty)\e^{-\gamma(\kappa)t},
    \qquad t\to\infty.
\]
Hence, for every fixed $s>0$,
\[
    \frac{\mathbb P_\mu(Z_{t+s}>0)}{\mathbb P_\mu(Z_t>0)}
    \longrightarrow \e^{-\gamma(\kappa)s},
    \qquad t\to\infty.
\]
Lemma \ref{lem_conditional_limit_qsd}(ii) identifies the decay rate as
$\gamma(\kappa)$. Corollary \ref{main_thm_1} ensures that this QSD is
unique at that decay rate.
\qed

\noindent{\it Proof of Theorem \ref{thm_complete_qsd_classification}(ii).~}
For $0<\kappa<\kappa_*$, Lemma \ref{lem_exact_power_construction} provides a
QSD with decay rate $\gamma(\kappa)$ and modified Laplace transform
$\Phi_\kappa$.

If $\kappa=\kappa_*$ and $\kappa_c\leq1$, the point-mass Yaglom theorem
\cite[Theorem~1.1]{LZZ} supplies the Yaglom limit.  Lemma
\ref{lem_conditional_limit_qsd}(i) makes it a QSD, while
\cite[Propositions~2.9 and~2.14]{LZZ} identify its modified Laplace transform
with $\Phi_{\kappa_c}$ in the weakly and intermediately subcritical cases,
respectively.

If $\kappa=\kappa_*=1<\kappa_c$, the strongly subcritical Yaglom theorem
\cite[Theorem 1.1]{LZZ} supplies the point-mass conditional limit.  Lemma
\ref{lem_conditional_limit_qsd}(i) makes it a QSD.  The explicit transform in
\cite[Theorem 1.1]{LZZ} is the displayed Kummer member $\Phi_1$.

In either case, let $\rho$ be the decay rate of the QSD just obtained.
Lemma \ref{lem_generator_laplace} shows that
$\Phi_{\kappa_*}$ solves \eqref{master_ode} with spectral parameter $\rho$, whereas the
Kummer equation shows that it solves the same equation with spectral parameter
$\gamma(\kappa_*)$.  Subtracting the two equations gives
\[
    \bigl(\rho-\gamma(\kappa_*)\bigr)\Phi_{\kappa_*}(\lambda)=0.
\]
Since $\Phi_{\kappa_*}$ is non-zero, $\rho=\gamma(\kappa_*)$.

Assertion (i) shows that any QSD with decay rate $\gamma(\kappa)$ has
modified Laplace transform $\Phi_\kappa$. Uniqueness of Laplace transforms
therefore makes the QSD constructed above the unique QSD with that decay rate.
This proves assertion (ii).
\qed

{
\section{Proof of the complete domains of attraction}\label{sec_complete_domains_attraction}
\setcounter{equation}{0}

{We now prove Theorem \ref{thm_complete_domains_attraction}, using the
QSD characterization in Theorem \ref{thm_complete_qsd_classification}.
The proof has three steps. Subsection \ref{sec_exponential_initial_masses}
establishes comparison and asymptotic estimates for deterministic initial
masses. Subsection \ref{sec_general_initial_laws} uses these estimates to
characterize conditional convergence for general initial laws in terms of
the logarithmic survival rate. Subsection \ref{sec_tail_criteria} then proves
that, for interior QSDs, the survival-rate condition is equivalent to a
logarithmic tail condition. Together with the endpoint moment criterion,
this completes the proof.}

\subsection{Exponentially growing initial masses}\label{sec_exponential_initial_masses}

{We establish two {results} for the treatment of general initial laws:
monotonicity of the conditional transform in the initial mass
(Lemma \ref{full_lem_initial_mass_monotonicity}), and survival asymptotics
and conditional limits for exponentially growing initial masses
(Lemma \ref{full_lem_moving_initial_mass}).}

{Recall that $V_t(\lambda)$ in \eqref{Vt_environment_definition}
is the random Laplace exponent given the environment. For an initial mass
$z>0$, put $f_z(r):=1-\e^{-zr}$ and write
\[
 A_t(z,\lambda):=\mbb E[f_z(V_t(\lambda))]
 =\mbb E_z[1-\e^{-\lambda Z_t}].
\]
Its value at $\lambda=\infty$ is the survival probability
$A_t(z,\infty)=\mbb P_z(Z_t>0)$. Dividing by this probability gives
the modified Laplace transform of the conditional law from $z$:
\[
 R_t(z,\lambda):=\frac{A_t(z,\lambda)}{A_t(z,\infty)}
 =\mbb E_z[1-\e^{-\lambda Z_t}\mid Z_t>0].
\]}

\bglemma
\label{full_lem_initial_mass_monotonicity}
For every fixed \(t,\lambda>0\), the function
\(z\mapsto R_t(z,\lambda)\) is nondecreasing.
\edlemma
\proof
Fix $0<z_1<z_2$. {By monotone convergence, for every $z>0$,
\[
 A_t(z,\lambda)=\mbb E_z[1-\e^{-\lambda Z_t}]
 \uparrow\mbb P_z(Z_t>0)=A_t(z,\infty)
 \qquad\text{as }\lambda\to\infty.
\]
Therefore, to prove the lemma, it suffices to show that, for every
$0<\lambda_1<\lambda_2<\infty$,
\begin{equation}\label{full_eq_cross_product}
 A_t(z_2,\lambda_1)A_t(z_1,\lambda_2)
 \geq A_t(z_1,\lambda_1)A_t(z_2,\lambda_2).
\end{equation}
Indeed, setting $\lambda_1=\lambda$, letting $\lambda_2\to\infty$, and
dividing by $A_t(z_1,\infty)A_t(z_2,\infty)>0$ gives
$R_t(z_2,\lambda)\geq R_t(z_1,\lambda)$.}

{Let $B$ be a standard Brownian motion and, for $\lambda>0$,
define $(Q_s^\lambda)_{s\geq0}$ by}
\begin{equation}\label{full_eq_auxiliary_solution}
 Q_s^\lambda
 \mathrel{:=}\e^{\sigma B_s+{\bf m}s}
  \left(\lambda^{-\beta}+c\beta\int_0^s
        \e^{\beta(\sigma B_r+{\bf m}r)}\,\d r\right)^{-1/\beta},
 \qquad s\geq0.
\end{equation}
{By \eqref{full_eq_physical_representation}, for each fixed $t>0$,
\[
 Q_t^\lambda\stackrel{\d}{=}V_t(\lambda),\qquad
 A_t(z,\lambda)=\mbb E[f_z(Q_t^\lambda)].
\]}
It\^o's formula, with $\alpha={\bf m}+\sigma^2/2$, gives
\begin{equation}\label{full_eq_auxiliary_sde}
 \d Q_s^\lambda=\sigma Q_s^\lambda\,\d B_s
       +\bigl(\alpha Q_s^\lambda-c(Q_s^\lambda)^{1+\beta}\bigr)\,\d s,
 \qquad Q_0^\lambda=\lambda.
\end{equation}
By \eqref{full_eq_auxiliary_solution}, the solution neither hits zero nor explodes in finite time.
The coefficients are locally Lipschitz on $(0,\infty)$, so this is the
unique global strong solution.

The ratio $f_{z_2}(r)/f_{z_1}(r)$ is
nonincreasing in \(r>0\), because
\[
 \frac{\d}{\d r}\log\frac{f_{z_2}(r)}{f_{z_1}(r)}
 =\frac{z_2}{\e^{z_2r}-1}-\frac{z_1}{\e^{z_1r}-1}\le0.
\]
{Define the bounded antisymmetric function}
\[
 {G(u,v):=f_{z_2}(u)f_{z_1}(v)-f_{z_1}(u)f_{z_2}(v)}
\]
{The preceding monotonicity gives $G(u,v)\geq0$ whenever $u<v$.}

{Take an independent coupling
$(Q_s^{(1),\lambda_1},Q_s^{(2),\lambda_2})_{s\geq0}$
of two solutions of \eqref{full_eq_auxiliary_sde}, with initial values
$Q_0^{(i),\lambda_i}=\lambda_i$, $i=1,2$.}
{By independence,
\begin{equation}\label{full_eq_coupling_cross_product}
\begin{aligned}
 \mbb E[G(Q_t^{(1),\lambda_1},Q_t^{(2),\lambda_2})]
 &={}
 A_t(z_2,\lambda_1)A_t(z_1,\lambda_2)\\
 &\quad-A_t(z_1,\lambda_1)A_t(z_2,\lambda_2).
\end{aligned}
\end{equation}
Thus \eqref{full_eq_cross_product} is equivalent to
$\mbb E[G(Q_t^{(1),\lambda_1},Q_t^{(2),\lambda_2})]\geq0$.}

Let $(\mathcal G_s)_{s\geq0}$ be their usual joint natural filtration,
and let $T$ be their first meeting time, with $T=\infty$ if they never
meet.

{On $\{T\leq t\}$, put
$q:=Q_T^{(1),\lambda_1}=Q_T^{(2),\lambda_2}$ and $s:=t-T$.
The strong Markov property of the pair, together with
$\{T\leq t\}\in\mathcal G_T$, gives
\begin{equation}\label{full_eq_coupling_markov}
\begin{aligned}
 &\mbb E\!\left[
 \mathbf1_{\{T\leq t\}}f_{z_j}(Q_t^{(i),\lambda_i})
 \mid\mathcal G_T\right]\\
 &\quad=\mathbf1_{\{T\leq t\}}\mbb E[f_{z_j}(Q_s^q)]
 =\mathbf1_{\{T\leq t\}}A_s(z_j,q),\qquad i,j\in\{1,2\}.
\end{aligned}
\end{equation}
Using \eqref{full_eq_coupling_markov} and the conditional independence
of the two copies given $\mathcal G_T$ on $\{T\leq t\}$, we obtain
\begin{equation}\label{full_eq_coupling_after_meeting}
\begin{aligned}
 &\mbb E\!\left[
     \mathbf1_{\{T\leq t\}}G(Q_t^{(1),\lambda_1},Q_t^{(2),\lambda_2})
     \mid\mathcal G_T
 \right]\\
 &\quad=\mathbf1_{\{T\leq t\}}\Bigl(
         \mbb E\!\left[f_{z_2}(Q_t^{(1),\lambda_1})\mid\mathcal G_T\right]
         \mbb E\!\left[f_{z_1}(Q_t^{(2),\lambda_2})\mid\mathcal G_T\right]\\
 &\qquad\qquad-\mbb E\!\left[f_{z_1}(Q_t^{(1),\lambda_1})\mid\mathcal G_T\right]
          \mbb E\!\left[f_{z_2}(Q_t^{(2),\lambda_2})\mid\mathcal G_T\right]\Bigr)\\
 &\quad=\mathbf1_{\{T\leq t\}}
 \bigl(A_s(z_2,q)A_s(z_1,q)-A_s(z_1,q)A_s(z_2,q)\bigr)=0.
\end{aligned}
\end{equation}
On $\{T>t\}$, the continuous paths cannot reverse their initial order
without meeting, so $Q_t^{(1),\lambda_1}<Q_t^{(2),\lambda_2}$.
The monotonicity of $f_{z_2}/f_{z_1}$ therefore gives
$G(Q_t^{(1),\lambda_1},Q_t^{(2),\lambda_2})\geq0$.
Taking expectations in \eqref{full_eq_coupling_after_meeting} now yields
\begin{align*}
 &\mbb E[G(Q_t^{(1),\lambda_1},Q_t^{(2),\lambda_2})]\\
 &\quad=\mbb E[\mathbf1_{\{T>t\}}G(Q_t^{(1),\lambda_1},Q_t^{(2),\lambda_2})]\geq0.
\end{align*}}

{By \eqref{full_eq_coupling_cross_product}, this proves
\eqref{full_eq_cross_product} and hence the assertion.}
\qed

{For the asymptotic estimates below, put
\[
 J(v):=\frac{(-{\bf m}-v)^2}{2\sigma^2}.
\]}
{
We first record the Brownian bridge estimates needed below.

\bglemma[Brownian bridge estimates]\label{full_lem_bridge_domination}
Let $W$ be a standard Brownian motion, and fix $v>0$ and $0<\kappa<1$.
For $t>0$, $y\in\mbb R$, and $\lambda\in(0,\infty]$, set
\[
\begin{aligned}
 D_t^y(\lambda)
 &:=\lambda^{-\beta}+c\beta\int_0^t
 \exp\left\{\beta\left(
 \sigma W_s-\frac{\sigma s}{t}W_t+\frac{sy}{t}-vs
 \right)\right\}\,\d s,\\
 D_\infty(\lambda)
 &:=\lambda^{-\beta}+c\beta\int_0^\infty
        \e^{\beta(\sigma W_s-vs)}\,\d s,
\end{aligned}
\]
with $\infty^{-\beta}=0$. For every fixed $y\in\mbb R$ and
$\lambda\in(0,\infty]$,
\[
 D_t^y(\lambda)\longrightarrow D_\infty(\lambda)\in(0,\infty)
 \qquad\text{almost surely as }t\to\infty.
\]
Moreover, there exist $t_0\geq1$ and a nonnegative function
$g\in L^1(\mbb R)$ such that
\[
 0\leq \e^{-\kappa y}
 \mbb E\!\left[f_1\!\left(\e^yD_t^y(\lambda)^{-1/\beta}\right)\right]
 \leq g(y),
 \qquad t\geq t_0,\quad y\in\mbb R,\quad\lambda\in(0,\infty].
\]
\edlemma

\proof
Fix $y\in\mbb R$ and work on the event of probability one where
$W_r/r\to0$ as $r\to\infty$. Set
\[
 h_t(s):=\mathbf1_{[0,t]}(s)
 \exp\left\{\beta\left(\sigma W_s-vs+
                  \frac{s}{t}(y-\sigma W_t)\right)\right\},
 \qquad s\geq0.
\]
For each fixed $s\geq0$,
\[
 \frac{s}{t}(y-\sigma W_t)\longrightarrow0,
 \qquad
 h_t(s)\longrightarrow\e^{\beta(\sigma W_s-vs)}.
\]
For all sufficiently large $t$, $|y-\sigma W_t|/t\leq v/2$, so
\begin{equation}\label{full_eq_bridge_integrand_bound}
\begin{aligned}
 0\leq h_t(s)
 &\leq\exp\left\{\beta\left(\sigma W_s-vs+
                  s\frac{|y-\sigma W_t|}{t}\right)\right\}\\
 &\leq\e^{\beta(\sigma W_s-vs/2)},\qquad s\geq0.
\end{aligned}
\end{equation}
Choose $S<\infty$ such that $\sigma W_s\leq vs/4$ for $s\geq S$.
Continuity of $W$ on $[0,S]$ gives
\[
\begin{aligned}
 \int_0^\infty\e^{\beta(\sigma W_s-vs/2)}\,\d s
 &\leq\int_0^S\e^{\beta(\sigma W_s-vs/2)}\,\d s
       +\int_S^\infty\e^{-\beta vs/4}\,\d s\\
 &=\int_0^S\e^{\beta(\sigma W_s-vs/2)}\,\d s
       +\frac4{\beta v}\e^{-\beta vS/4}<\infty.
\end{aligned}
\]
Dominated convergence, using \eqref{full_eq_bridge_integrand_bound}, now yields
\[
\begin{aligned}
 D_t^y(\lambda)
 &=\lambda^{-\beta}+c\beta\int_0^\infty h_t(s)\,\d s\\
 &\longrightarrow\lambda^{-\beta}+c\beta\int_0^\infty
                  \e^{\beta(\sigma W_s-vs)}\,\d s
 =D_\infty(\lambda)\in(0,\infty),
\end{aligned}
\]
which proves the first assertion.

For $t\geq1$, restricting the integral defining $D_t^y(\lambda)$ to
$[0,1]$ gives
\[
 D_t^y(\lambda)^{-1/\beta}
 \leq (c\beta)^{-1/\beta}\e^v
 \exp\left\{\sigma\sup_{r\leq1}|W_r|
              +\frac{\sigma|W_t|}{t}+\frac{|y|}{t}\right\}.
\]
Cauchy--Schwarz, the Gaussian tails of $\sup_{r\leq1}|W_r|$, and
$W_t/\sqrt t\sim N(0,1)$ yield a constant $C<\infty$, independent of
$t,y,\lambda$, such that
\[
 \mbb E[D_t^y(\lambda)^{-1/\beta}]
 \leq C\e^{|y|/t},\qquad t\geq1.
\]
Using $0\leq f_1(u)\leq\min\{1,u\}$, we obtain
\[
 \e^{-\kappa y}
 \mbb E\!\left[f_1\!\left(\e^yD_t^y(\lambda)^{-1/\beta}\right)\right]
 \leq \e^{-\kappa y}\min\{1,C\e^{y+|y|/t}\}.
\]
For $t\geq t_0:=2/(1-\kappa)$, the right-hand side is bounded by
\[
 g(y):=
 \begin{cases}
  \e^{-\kappa y},&y\geq0,\\[1mm]
  C\e^{(1-\kappa)y/2},&y<0,
 \end{cases}
 \qquad
 \int_{\mbb R}g(y)\,\d y
 =\frac1\kappa+\frac{2C}{1-\kappa}<\infty.
\]
This proves the second assertion.
\qed
}

We now compute the asymptotics for initial mass \(\e^{vt}\).  This is the
scale that will later be matched with the tail of the initial distribution.  A
tail with logarithmic exponent $\kappa$ assigns to masses of order $\e^{vt}$
a cost of order $\e^{-\kappa vt}$ on the logarithmic scale, while survival
from such a moving initial mass has its own exponential cost.  The balance of
these two contributions is what determines the relevant initial scale in the
interior attraction problem.

\bglemma
\label{full_lem_moving_initial_mass}
{Let \(v\) be fixed with
\(-{\bf m}-\kappa_*\sigma^2<v<-{\bf m}\), and set
\[
 \kappa:=\frac{-{\bf m}-v}{\sigma^2}\in(0,\kappa_*).
\]}
For every
\(\lambda\in(0,\infty]\),
\begin{equation}\label{full_eq_moving_asymptotic}
 A_t(\e^{vt},\lambda)
 \sim\frac{\Gamma(1-\kappa)}{\kappa\sigma\sqrt{2\pi}}
       F_{\kappa/\beta}(\lambda)t^{-1/2}\e^{-J(v)t},
 \qquad t\to\infty.
\end{equation}
In particular, for every fixed \(\lambda>0\),
\begin{equation}\label{full_eq_moving_transform_limit}
 {\lim_{t\to\infty}R_t(\e^{vt},\lambda)
 =\Phi_{(-{\bf m}-v)/\sigma^2}(\lambda).}
\end{equation}
\edlemma
\proof
By the definition of $A_t$, the identity $f_z(r)=f_1(zr)$, and
\eqref{full_eq_physical_representation},
\begin{equation}\label{full_eq_moving_Y_representation}
\begin{aligned}
 A_t(\e^{vt},\lambda)
 &=\mbb E[f_{\e^{vt}}(V_t(\lambda))]\\
 &=\mbb E[f_1(\e^{vt}V_t(\lambda))]\\
 &=\mbb E\left[
 f_1\left(
   \e^{\sigma B_t+({\bf m}+v)t}
   \left(\lambda^{-\beta}+c\beta\int_0^t
      \e^{\beta(\sigma B_s+{\bf m}s)}\,\d s
   \right)^{-1/\beta}
 \right)\right]\\
 &=\mbb E\left[
 f_1\left(
   \e^{Y_t}
   \left(\lambda^{-\beta}+c\beta\int_0^t
      \e^{\beta(Y_s-vs)}\,\d s
   \right)^{-1/\beta}
 \right)\right].
\end{aligned}
\end{equation}
Here
\[
 Y_s=\sigma B_s+({\bf m}+v)s,\qquad s\geq0.
\]
Since ${\bf m}+v=-\kappa\sigma^2$, the density of $Y_t$ is
\[
 p_t(y):=
 \frac{\e^{-J(v)t}}{\sigma\sqrt{2\pi t}}
 \exp\left\{-\frac{y^2}{2\sigma^2t}-\kappa y\right\},
 \qquad y\in\mbb R.
\]
Conditionally on $Y_t=y$, the path $(Y_s)_{0\leq s\leq t}$ has the
same law as
\[
 \left(
 \sigma\left(W_s-\frac{s}{t}W_t\right)+\frac{sy}{t}
 \right)_{0\leq s\leq t},
\]
{where $W$ is a standard Brownian motion}; see
\cite[Section~3.4]{PY18} for the Brownian bridge representation.
 Consequently, conditionally on
$Y_t=y$,
\[
 Y_s-vs
 \overset{\mathrm d}{=}
 \sigma W_s-\frac{\sigma s}{t}W_t+\frac{sy}{t}-vs.
\]
{With $D_t^y(\lambda)$ as in Lemma \ref{full_lem_bridge_domination},
we therefore have, conditionally on $Y_t=y$,}
\[
 \lambda^{-\beta}
 +c\beta\int_0^t\e^{\beta(Y_s-vs)}\,\d s
 \ \overset{\mathrm d}{=}\
 D_t^y(\lambda).
\]
Applying the law of total expectation to
\eqref{full_eq_moving_Y_representation} and then using the density $p_t$
of $Y_t$ gives
\begin{equation}\label{full_eq_bridge_integral}
\begin{aligned}
 &\sqrt t\,\e^{J(v)t}A_t(\e^{vt},\lambda)\\
 &\quad=\sqrt t\,\e^{J(v)t}\int_{\mbb R}
 \mbb E\left[
 f_1\bigl(\e^yD_t^y(\lambda)^{-1/\beta}\bigr)
 \right]p_t(y)\,\d y\\
 &\quad=\frac1{\sigma\sqrt{2\pi}}\int_{\mbb R}
 \e^{-y^2/(2\sigma^2t)}\e^{-\kappa y}
 \mbb E\left[
 f_1\bigl(\e^yD_t^y(\lambda)^{-1/\beta}\bigr)
 \right]\,\d y.
\end{aligned}
\end{equation}
{Lemma \ref{full_lem_bridge_domination} gives, for each fixed $y$,
\[
 D_t^y(\lambda)\longrightarrow D_\infty(\lambda)
 \qquad\text{almost surely as }t\to\infty,
\]
together with an integrable bound for the integrand in
\eqref{full_eq_bridge_integral}, independent of all sufficiently large $t$.
Bounded convergence inside the expectation and dominated convergence in $y$
therefore give
\begin{equation}\label{full_eq_moving_limit_integral}
\begin{aligned}
&\lim_{t\to\infty}
 \sqrt t\,\e^{J(v)t}A_t(\e^{vt},\lambda)\\
&\quad=
 \frac1{\sigma\sqrt{2\pi}}
 \int_{\mbb R}\e^{-\kappa y}
 \mbb E\!\left[
 1-\e^{-\e^yD_\infty(\lambda)^{-1/\beta}}
 \right]\,\d y\\
&\quad=
 \frac{\Gamma(1-\kappa)}{\kappa\sigma\sqrt{2\pi}}
 \mbb E[D_\infty(\lambda)^{-\kappa/\beta}]\\
&\quad=
 \frac{\Gamma(1-\kappa)}{\kappa\sigma\sqrt{2\pi}}
 F_{\kappa/\beta}(\lambda).
\end{aligned}
\end{equation}
For the second equality, we use Tonelli's theorem and the substitution
$u=\e^yD_\infty(\lambda)^{-1/\beta}$, since $D_\infty(\lambda)$ does not
depend on $y$. This proves \eqref{full_eq_moving_asymptotic}.

Dividing \eqref{full_eq_moving_asymptotic} by its value at $\lambda=\infty$
and using Lemma \ref{lem_fractional_moment_strict} gives
\[
 \lim_{t\to\infty}R_t(\e^{vt},\lambda)
 =\frac{F_{\kappa/\beta}(\lambda)}{F_{\kappa/\beta}(\infty)}
 =\Phi_{(-{\bf m}-v)/\sigma^2}(\lambda),
\]
which proves \eqref{full_eq_moving_transform_limit}.}
\qed

\subsection{General initial laws conditioned on survival}\label{sec_general_initial_laws}

{We now establish convergence of the conditional law from general initial
distributions to the interior QSDs under a survival-rate condition
(Proposition \ref{full_prop_interior_concentration}), and to the endpoint QSD
under a moment condition (Proposition \ref{full_prop_endpoint_attraction}).
These results lead to a characterization of conditional convergence in
terms of the logarithmic survival rate
(Corollary \ref{full_cor_survival_criterion}).}

For an initial probability \(\mu\) on \((0,\infty)\), write
\[
 \mu_t=\mbb P_\mu(Z_t\in\cdot\mid Z_t>0).
\]

{
We use the convention $-\log0=\infty$.

\bglemma\label{full_lem_tail_moment}
For every probability measure $\mu$ on $(0,\infty)$,
\[
 \sup\left\{a\geq0:\int_0^\infty z^a\mu(\d z)<\infty\right\}
 =\liminf_{x\to\infty}\frac{-\log\mu([x,\infty))}{\log x}.
\]
\edlemma
\proof
For
\[
 0<a<b<\liminf_{x\to\infty}\frac{-\log\mu([x,\infty))}{\log x},
\]
choose $x_0\geq1$ such that $\mu([x,\infty))\leq x^{-b}$ for all
$x\geq x_0$. Tonelli's theorem gives
\[
\begin{aligned}
 \int_0^\infty z^a\mu(\d z)
 &=\int_0^\infty\left(\int_0^z a x^{a-1}\,\d x\right)\mu(\d z)\\
 &=a\int_0^\infty x^{a-1}\mu([x,\infty))\,\d x\\
 &\leq x_0^a+a\int_{x_0}^\infty x^{a-b-1}\,\d x<\infty.
\end{aligned}
\]
Thus every positive moment of order below the lower logarithmic tail
exponent is finite.

Conversely, if $a>0$ and $\int_0^\infty z^a\mu(\d z)<\infty$,
Markov's inequality gives
\[
 \mu([x,\infty))\leq x^{-a}\int_0^\infty z^a\mu(\d z),
 \qquad
 \liminf_{x\to\infty}\frac{-\log\mu([x,\infty))}{\log x}\geq a.
\]
Taking the supremum over finite moment orders proves the reverse
inequality; if no positive moment is finite, it follows from
$\mu([x,\infty))\leq1$.
\qed}

{The following bound allows us to use different moment orders when
averaging over the initial law.

\bglemma[Uniform survival bound]\label{full_lem_uniform_survival_bound}
For every fixed $0<\theta<\kappa_*$, there is a constant
$C_\theta<\infty$ such that, for all sufficiently large $t$,
\begin{equation}\label{full_eq_uniform_survival_moment_bound}
 \mbb P_z(Z_t>0)
 \leq C_\theta z^\theta\e^{-\gamma(\theta)t},\qquad z>0.
\end{equation}
Both the constant and the lower threshold for $t$ may depend on $\theta$,
but are independent of $z$.
\edlemma
\proof
Since $0<\theta<\kappa_*\leq1$,
\[
 1-\e^{-u}\leq\min\{u,1\}\leq u^\theta,
 \qquad u\geq0.
\]
The Laplace formula and monotone convergence therefore give
\[
\begin{aligned}
 \mbb P_z(Z_t>0)
 &=\lim_{\lambda\to\infty}\mbb E_z[1-\e^{-\lambda Z_t}]\\
 &=\mbb E[1-\e^{-zV_t(\infty)}]
 \leq z^\theta\mbb E[V_t(\infty)^\theta].
\end{aligned}
\]
To estimate the remaining moment, note that
$\theta<\kappa_*\leq\kappa_c=-{\bf m}/\sigma^2$, so
\[
 0<\frac{\theta}{\beta}
   <\frac{-{\bf m}}{\beta\sigma^2}=\frac{\eta}{2}.
\]
By \eqref{Vt_explicit_representation} and the definition of
$\mathcal M_{p,t}$ in Lemma \ref{lem_fractional_moment_strict},
\[
 \mbb E[V_t(\infty)^\theta]
 =\mbb E\!\left[(\gamma_0A_{\tau_t}^{(\eta)})^{-\theta/\beta}\right]
 =\mathcal M_{\theta/\beta,t}(\infty).
\]
Applying that lemma with $p=\theta/\beta$ and $\lambda=\infty$ yields
\[
 \e^{\gamma(\theta)t}\mbb E[V_t(\infty)^\theta]
 \longrightarrow F_{\theta/\beta}(\infty)\in(0,\infty),
 \qquad t\to\infty.
\]
Thus, with $C_\theta:=2F_{\theta/\beta}(\infty)$, for all sufficiently
large $t$,
\[
 \mbb E[V_t(\infty)^\theta]
 \leq C_\theta\e^{-\gamma(\theta)t}.
\]
Neither this constant nor the time threshold depends on $z$.
Substituting this moment bound into the preceding estimate proves
\eqref{full_eq_uniform_survival_moment_bound}.
\qed

\bglemma\label{full_lem_conditional_survival_rate}
Let $\mu$ be a probability measure on $(0,\infty)$.
If $\mu_t\Longrightarrow\nu_\kappa$ for some
$0<\kappa\leq\kappa_*$, then
\begin{equation}\label{full_eq_rate_from_conditional_limit}
 \lim_{t\to\infty}-\frac1t\log \mbb P_\mu(Z_t>0)=\gamma(\kappa).
\end{equation}
In particular, for every fixed $z>0$,
\begin{equation}\label{full_eq_fixed_survival_rate}
 \lim_{t\to\infty}-\frac1t\log\mbb P_z(Z_t>0)=\gamma(\kappa_*).
\end{equation}
\edlemma
\proof
Suppose $\mu_t\Longrightarrow\nu_\kappa$.
For each $h>0$, the function $z\mapsto\mbb P_z(Z_h>0)$ is bounded
and continuous. The Markov property gives
\[
 \frac{\mbb P_\mu(Z_{t+h}>0)}{\mbb P_\mu(Z_t>0)}
 =\int_0^\infty\mbb P_z(Z_h>0)\,\mu_t(\d z)
 \longrightarrow
 \int_0^\infty\mbb P_z(Z_h>0)\,\nu_\kappa(\d z)
 =\e^{-\gamma(\kappa)h}.
\]
Take \(h=1\). The preceding limit gives
\[
 \frac{\mbb P_\mu(Z_{n+1}>0)}{\mbb P_\mu(Z_n>0)}
 \longrightarrow\e^{-\gamma(\kappa)}.
\]
Since \(\mbb P_\mu(Z_0>0)=1\),
\[
 -\frac1n\log \mbb P_\mu(Z_n>0)
 =
 -\frac1n\sum_{j=0}^{n-1}
       \log\frac{\mbb P_\mu(Z_{j+1}>0)}{\mbb P_\mu(Z_j>0)}
 \longrightarrow
 -\log\e^{-\gamma(\kappa)}
 =
 \gamma(\kappa).
\]
Monotonicity of $t\mapsto\mbb P_\mu(Z_t>0)$ extends this limit
from integer to arbitrary times, proving
\eqref{full_eq_rate_from_conditional_limit}.
The conditional law from each $z>0$ converges to $\nu_{\kappa_*}$ by
\cite[Theorem~1.1]{LZZ}. Thus \eqref{full_eq_rate_from_conditional_limit},
with $\mu=\delta_z$ and $\kappa=\kappa_*$, gives
\eqref{full_eq_fixed_survival_rate}.
\qed

For any initial probability measure $\mu$, choose $z_0>0$ with
$\mu([z_0,\infty))>0$. Monotonicity gives
\[
 \mbb P_\mu(Z_t>0)
 \geq\mu([z_0,\infty))\mbb P_{z_0}(Z_t>0).
\]
Consequently, \eqref{full_eq_fixed_survival_rate} implies the general bound
\begin{equation}\label{full_eq_universal_survival_rate}
 \limsup_{t\to\infty}-\frac1t\log\mbb P_\mu(Z_t>0)
 \leq\gamma(\kappa_*).
\end{equation}

\bglemma\label{full_lem_survival_rates}
Let $\mu$ be an initial probability measure on $(0,\infty)$.
\begin{itemize}
\item[\rm(i)] If
\[
 \int_0^\infty z^a\mu(\d z)<\infty,
 \qquad\text{for every }a\in(0,\kappa_*),
\]
then
\[
 \lim_{t\to\infty}-\frac1t\log\mbb P_\mu(Z_t>0)=\gamma(\kappa_*).
\]
\item[\rm(ii)] If, for some $0<\kappa\leq\kappa_*$,
\[
 \lim_{t\to\infty}-\frac1t\log\mbb P_\mu(Z_t>0)=\gamma(\kappa),
\]
then
\[
 \int_0^\infty z^a\mu(\d z)<\infty,
 \qquad\text{for every }a\in(0,\kappa).
\]
\end{itemize}
\edlemma
\proof
For (i), fix $0<\theta<\kappa_*$. By
\eqref{full_eq_uniform_survival_moment_bound} and the moment assumption,
for all sufficiently large $t$,
\[
\begin{aligned}
 \mbb P_\mu(Z_t>0)
 &=\int_0^\infty\mbb P_z(Z_t>0)\mu(\d z)\\
 &\leq C_\theta\left(\int_0^\infty z^\theta\mu(\d z)\right)
       \e^{-\gamma(\theta)t}<\infty.
\end{aligned}
\]
Taking logarithms and using \eqref{full_eq_universal_survival_rate} gives
\[
 \gamma(\theta)
 \leq\liminf_{t\to\infty}-\frac1t\log\mbb P_\mu(Z_t>0)
 \leq\limsup_{t\to\infty}-\frac1t\log\mbb P_\mu(Z_t>0)
 \leq\gamma(\kappa_*).
\]
Letting $\theta\uparrow\kappa_*$ proves (i).

For (ii), fix $0<a<\kappa$. By monotonicity and Lemma
\ref{full_lem_moving_initial_mass},
\[
\begin{aligned}
 \mu([\e^{(-{\bf m}-a\sigma^2)t},\infty))
 &\leq\frac{\mbb P_\mu(Z_t>0)}
              {\mbb P_{\e^{(-{\bf m}-a\sigma^2)t}}(Z_t>0)}\\
 &=\exp\!\left\{-\bigl[\gamma(\kappa)
                -J(-{\bf m}-a\sigma^2)\bigr]t+o(t)\right\}.
\end{aligned}
\]
Since $a(-{\bf m}-a\sigma^2)+J(-{\bf m}-a\sigma^2)=\gamma(a)$,
setting $x=\e^{(-{\bf m}-a\sigma^2)t}$ gives
\[
 \mu([x,\infty))
 \leq x^{-a-\frac{\gamma(\kappa)-\gamma(a)}{-{\bf m}-a\sigma^2}+o(1)},
 \qquad x\to\infty.
\]
Here $\gamma(\kappa)>\gamma(a)$ and $-{\bf m}-a\sigma^2>0$.
Thus the tail integral
formula in the proof of Lemma \ref{full_lem_tail_moment} gives
$\int_0^\infty z^a\mu(\d z)<\infty$.
\qed}

To prove conditional convergence, we locate the initial masses that
contribute to survival.  The contribution of an initial mass $z$ is weighted
by both its probability under $\mu$ and its probability of survival up to time
$t$.  Thus conditioning on survival replaces the original initial law by the
following survival-biased law. Define
\[
 \pi_t(\d z)
 :=\mbb P_\mu(Z_0\in\d z\mid Z_t>0)
 =\frac{\mbb P_z(Z_t>0)\mu(\d z)}{\mbb P_\mu(Z_t>0)}.
\]
Thus $\pi_t$ is the survival-biased distribution of the initial mass,
whereas $\mu_t$ is {the conditional law from $\mu$ at time $t$}.

The annealed ratio \eqref{annealed_ratio} then reads
\begin{equation}\label{full_eq_conditional_mixture}
 \mu_t(f_\lambda)=\int_0^\infty R_t(z,\lambda)\pi_t(\d z),
 \qquad \lambda>0.
\end{equation}

We will also use continuity of \(\kappa\mapsto\Phi_\kappa(\lambda)\)
on \((0,\kappa_*]\), for each fixed \(\lambda>0\). To see this,
put \(z=2c\lambda^\beta/(\beta\sigma^2)>0\) in its Kummer expression
and use \cite[(13.4.4)]{NISTDLMF}:
\[
 U(a,b,z)=\frac1{\Gamma(a)}\int_0^\infty
       \e^{-zu}u^{a-1}(1+u)^{b-a-1}\,\d u,\qquad a>0.
\]
For \(\kappa\) in a compact subinterval of \((0,\kappa_*]\), the
integrand has an integrable bound at both endpoints. Dominated convergence
therefore also gives continuity at \(\kappa_*\).

{The following proposition proves concentration of the initial mass
and conditional convergence under a survival-rate assumption.}

\bgproposition[{Interior concentration and attraction}]
\label{full_prop_interior_concentration}
{Let $0<\kappa<\kappa_*$ and suppose that}
\begin{equation}\label{full_eq_interior_survival_rate}
 \lim_{t\to\infty}-\frac1t\log \mbb P_\mu(Z_t>0)=\gamma(\kappa).
\end{equation}
{For every sufficiently small $\varepsilon>0$,}
\begin{equation}\label{full_eq_interior_concentration}
 \lim_{t\to\infty}
 \pi_t\left(\left[\e^{((-{\bf m}-\kappa\sigma^2)-\varepsilon)t},
                         \e^{((-{\bf m}-\kappa\sigma^2)+\varepsilon)t}\right]\right)
 =1.
\end{equation}
{Moreover, $\mu_t\Longrightarrow\nu_\kappa$.}
\edproposition
\proof
{By \eqref{full_eq_interior_survival_rate} and Lemma
\ref{full_lem_survival_rates}(ii),
\begin{equation}\label{full_eq_interior_finite_moments}
 \int_0^\infty z^a\mu(\d z)<\infty,
 \qquad\text{for every }a\in(0,\kappa).
\end{equation}
Recall from Lemma \ref{full_lem_uniform_survival_bound} that, for every
$0<\theta<\kappa_*$ and all sufficiently large $t$,
\[
 \mbb P_z(Z_t>0)\leq C_\theta z^\theta\e^{-\gamma(\theta)t},
 \qquad z>0,
\]
where $C_\theta$ is independent of $z$ and $t$.
{Fix a sufficiently small $\varepsilon>0$.
Choose $\kappa-\varepsilon/\sigma^2<a<\kappa$ sufficiently close to
$\kappa$ that}
\begin{equation}\label{full_eq_interior_rate_gaps}
\begin{aligned}
 {\delta_{\mathrm R}}&:=\frac{\varepsilon^2}{2\sigma^2}
       -(\kappa-a)(-{\bf m}-\kappa\sigma^2+\varepsilon){{}>0},\\
 {\delta_{\mathrm L}}&:=\frac{\varepsilon^2}{2\sigma^2}
       -(\kappa-a)(-{\bf m}-\kappa\sigma^2-\varepsilon){{}\geq\delta_{\mathrm R}}.
\end{aligned}
\end{equation}
{Here the subscripts $\mathrm R$ and $\mathrm L$ refer to the right
and left tails, respectively.}
Apply \eqref{full_eq_uniform_survival_moment_bound} with
$\theta=\kappa-\varepsilon/\sigma^2$.
Since $\kappa-\varepsilon/\sigma^2-a<0$, we obtain
\begin{equation}\label{full_eq_interior_upper_survival_mass}
\begin{aligned}
 &\int_{z>\e^{(-{\bf m}-\kappa\sigma^2+\varepsilon)t}}
      \mbb P_z(Z_t>0)\mu(\d z)\\
 &\quad\leq C_{\kappa-\varepsilon/\sigma^2}
       \e^{-\gamma(\kappa-\varepsilon/\sigma^2)t}
       \int_{z>\e^{(-{\bf m}-\kappa\sigma^2+\varepsilon)t}}
          z^{\kappa-\varepsilon/\sigma^2}\mu(\d z)\\
 &\quad= C_{\kappa-\varepsilon/\sigma^2}
       \e^{-\gamma(\kappa-\varepsilon/\sigma^2)t}
       \int_{z>\e^{(-{\bf m}-\kappa\sigma^2+\varepsilon)t}}
          z^a z^{\kappa-\varepsilon/\sigma^2-a}\mu(\d z)\\
 &\quad\leq C_{\kappa-\varepsilon/\sigma^2}
       \left(\int_0^\infty z^a\mu(\d z)\right)
       \e^{-[\gamma(\kappa-\varepsilon/\sigma^2)
                +(a-\kappa+\varepsilon/\sigma^2)(-{\bf m}-\kappa\sigma^2+\varepsilon)]t}\\
 &\quad= C_{\kappa-\varepsilon/\sigma^2}
       \left(\int_0^\infty z^a\mu(\d z)\right)
       \e^{-[\gamma(\kappa)+{\delta_{\mathrm R}}]t}<\infty,
\end{aligned}
\end{equation}
where the last equality uses \eqref{full_eq_interior_rate_gaps} and the last
inequality follows from \eqref{full_eq_interior_finite_moments}.

Next apply \eqref{full_eq_uniform_survival_moment_bound} with
$\theta=\kappa+\varepsilon/\sigma^2$.
Since $\kappa+\varepsilon/\sigma^2-a>0$, we similarly have
\begin{equation}\label{full_eq_interior_lower_survival_mass}
\begin{aligned}
 &\int_{z<\e^{(-{\bf m}-\kappa\sigma^2-\varepsilon)t}}
      \mbb P_z(Z_t>0)\mu(\d z)\\
 &\quad\leq C_{\kappa+\varepsilon/\sigma^2}
       \e^{-\gamma(\kappa+\varepsilon/\sigma^2)t}
       \int_{z<\e^{(-{\bf m}-\kappa\sigma^2-\varepsilon)t}}
          z^{\kappa+\varepsilon/\sigma^2}\mu(\d z)\\
 &\quad= C_{\kappa+\varepsilon/\sigma^2}
       \e^{-\gamma(\kappa+\varepsilon/\sigma^2)t}
       \int_{z<\e^{(-{\bf m}-\kappa\sigma^2-\varepsilon)t}}
          z^a z^{\kappa+\varepsilon/\sigma^2-a}\mu(\d z)\\
 &\quad\leq C_{\kappa+\varepsilon/\sigma^2}
       \left(\int_0^\infty z^a\mu(\d z)\right)
       \e^{-[\gamma(\kappa+\varepsilon/\sigma^2)
                +(a-\kappa-\varepsilon/\sigma^2)(-{\bf m}-\kappa\sigma^2-\varepsilon)]t}\\
 &\quad= C_{\kappa+\varepsilon/\sigma^2}
       \left(\int_0^\infty z^a\mu(\d z)\right)
       \e^{-[\gamma(\kappa)+{\delta_{\mathrm L}}]t}<\infty,
\end{aligned}
\end{equation}
{where the last equality uses \eqref{full_eq_interior_rate_gaps} and the last
inequality follows from \eqref{full_eq_interior_finite_moments}.}

By \eqref{full_eq_interior_survival_rate}, for all sufficiently large $t$,
\begin{equation}\label{full_eq_interior_survival_lower_bound}
 \mbb P_\mu(Z_t>0)\geq\e^{-[\gamma(\kappa)+{\delta_{\mathrm R}}/2]t}.
\end{equation}
{Set}
\[
 I_t:=\left[\e^{(-{\bf m}-\kappa\sigma^2-\varepsilon)t},
             \e^{(-{\bf m}-\kappa\sigma^2+\varepsilon)t}\right],
 \qquad {I_t^c:=(0,\infty)\setminus I_t}.
\]
{Adding \eqref{full_eq_interior_upper_survival_mass} and
\eqref{full_eq_interior_lower_survival_mass}, then using
$\delta_{\mathrm L}\geq\delta_{\mathrm R}$ and
\eqref{full_eq_interior_survival_lower_bound}, gives
\begin{equation}\label{full_eq_interior_normalized_tails}
\begin{aligned}
 \pi_t(I_t^c)
 &=\frac{\displaystyle\int_{I_t^c}\mbb P_z(Z_t>0)\mu(\d z)}
       {\mbb P_\mu(Z_t>0)}\\
 &\leq\left(C_{\kappa-\varepsilon/\sigma^2}
            +C_{\kappa+\varepsilon/\sigma^2}\right)
       \left(\int_0^\infty z^a\mu(\d z)\right)
       \frac{\e^{-[\gamma(\kappa)+\delta_{\mathrm R}]t}}
            {\e^{-[\gamma(\kappa)+\delta_{\mathrm R}/2]t}}\\
 &=\left(C_{\kappa-\varepsilon/\sigma^2}
         +C_{\kappa+\varepsilon/\sigma^2}\right)
       \left(\int_0^\infty z^a\mu(\d z)\right)
       \e^{-\delta_{\mathrm R}t/2}
       \xrightarrow[t\to\infty]{}0.
\end{aligned}
\end{equation}
The limit uses \eqref{full_eq_interior_finite_moments} and $\delta_{\mathrm R}>0$.}
This proves \eqref{full_eq_interior_concentration}.

{Fix $\lambda>0$.}
Since $\pi_t$ is a probability measure, \eqref{full_eq_conditional_mixture}
gives
\[
\begin{aligned}
 \mu_t(f_\lambda)-\Phi_\kappa(\lambda)
 &=\int_0^\infty R_t(z,\lambda)\pi_t(\d z)
   -\Phi_\kappa(\lambda)\int_0^\infty\pi_t(\d z)\\
 &=\int_0^\infty
       [R_t(z,\lambda)-\Phi_\kappa(\lambda)]\pi_t(\d z).
\end{aligned}
\]
The monotonicity of $R_t(\cdot,\lambda)$ in Lemma
\ref{full_lem_initial_mass_monotonicity} gives
\[
 |R_t(z,\lambda)-\Phi_\kappa(\lambda)|
 \leq\max_{\pm}\bigl|
       R_t(\e^{(-{\bf m}-\kappa\sigma^2\pm\varepsilon)t},\lambda)
       -\Phi_\kappa(\lambda)\bigr|,
 \qquad z\in I_t.
\]
Since $R_t(z,\lambda),\Phi_\kappa(\lambda)\in[0,1]$, splitting the integral
yields
\begin{equation}\label{full_eq_conditional_transform_error}
\begin{aligned}
 \bigl|\mu_t(f_\lambda)-\Phi_\kappa(\lambda)\bigr|
 &\leq\int_{{I_t^c}}
       |R_t(z,\lambda)-\Phi_\kappa(\lambda)|\pi_t(\d z)\\
 &\quad+\int_{I_t}
       |R_t(z,\lambda)-\Phi_\kappa(\lambda)|\pi_t(\d z)\\
 &\leq\pi_t\bigl({I_t^c}\bigr)\\
 &\quad+\pi_t(I_t)\max_{\pm}\bigl|
       R_t(\e^{(-{\bf m}-\kappa\sigma^2\pm\varepsilon)t},\lambda)
       -\Phi_\kappa(\lambda)\bigr|\\
 &\leq\pi_t\bigl({I_t^c}\bigr)\\
 &\quad+\max_{\pm}\bigl|
       R_t(\e^{(-{\bf m}-\kappa\sigma^2\pm\varepsilon)t},\lambda)
       -\Phi_\kappa(\lambda)\bigr|.
\end{aligned}
\end{equation}
By {\eqref{full_eq_interior_normalized_tails}},
$\pi_t({I_t^c})\to0$. Our choice of $\varepsilon$ ensures that
\[
 -{\bf m}-\kappa_*\sigma^2
 <-{\bf m}-\kappa\sigma^2-\varepsilon
 <-{\bf m}-\kappa\sigma^2+\varepsilon
 <-{\bf m}.
\]
Applying \eqref{full_eq_moving_transform_limit} first with
$v=-{\bf m}-\kappa\sigma^2-\varepsilon$ and then with
$v=-{\bf m}-\kappa\sigma^2+\varepsilon$, we obtain
\begin{equation}\label{full_eq_interior_boundary_transform_limits}
\begin{aligned}
 \lim_{t\to\infty}
 R_t(\e^{(-{\bf m}-\kappa\sigma^2-\varepsilon)t},\lambda)
 &=\Phi_{\kappa+\varepsilon/\sigma^2}(\lambda),\\
 \lim_{t\to\infty}
 R_t(\e^{(-{\bf m}-\kappa\sigma^2+\varepsilon)t},\lambda)
 &=\Phi_{\kappa-\varepsilon/\sigma^2}(\lambda).
\end{aligned}
\end{equation}
Thus \eqref{full_eq_conditional_transform_error},
\eqref{full_eq_interior_boundary_transform_limits}, and continuity of
$\kappa\mapsto\Phi_\kappa(\lambda)$ imply
\[
 \limsup_{t\to\infty}
       |\mu_t(f_\lambda)-\Phi_\kappa(\lambda)|
 \leq\max_{\pm}
       |\Phi_{\kappa\mp\varepsilon/\sigma^2}(\lambda)
          -\Phi_\kappa(\lambda)|
 \xrightarrow[\varepsilon\downarrow0]{}0.
\]
The continuity theorem for Laplace transforms proves
$\mu_t\Longrightarrow\nu_\kappa$.
\qed}

{
\bgproposition\label{full_prop_endpoint_attraction}
If
\begin{equation}\label{full_eq_endpoint_moment_condition}
 \int_0^\infty z^a\mu(\d z)<\infty,
 \qquad\text{for every }a\in(0,\kappa_*),
\end{equation}
then $\mu_t\Longrightarrow\nu_{\kappa_*}$.
\edproposition
\proof
Fix \(-{\bf m}-\kappa_*\sigma^2<v<-{\bf m}\).
Choose \(0<\theta<a<\kappa_*\) such that
\[
 \gamma(\theta)+(a-\theta)v>\gamma(\kappa_*).
\]
Such a choice exists because
\[
 \gamma(\kappa_*)-\gamma(\theta)
 =(-{\bf m}-\kappa_*\sigma^2)(\kappa_*-\theta)
    +\frac12\sigma^2(\kappa_*-\theta)^2:
\]
first choose \(\theta\) close to \(\kappa_*\), then choose \(a\)
sufficiently close to \(\kappa_*\) with \(a>\theta\).
By \eqref{full_eq_uniform_survival_moment_bound} in Lemma
\ref{full_lem_uniform_survival_bound}, for all sufficiently large $t$,
\[
\begin{aligned}
 \int_{z>\e^{vt}}\mbb P_z(Z_t>0)\mu(\d z)
 &\leq C_\theta\e^{-\gamma(\theta)t}
       \int_{z>\e^{vt}}z^\theta\mu(\d z)\\
 &=C_\theta\e^{-\gamma(\theta)t}
       \int_{z>\e^{vt}}z^a z^{\theta-a}\mu(\d z)\\
 &\leq C_\theta\left(\int_0^\infty z^a\mu(\d z)\right)
       \e^{-[\gamma(\theta)+(a-\theta)v]t}<\infty,
\end{aligned}
\]
where the second inequality uses
$z^{\theta-a}\leq\e^{-(a-\theta)vt}$ for $z>\e^{vt}$, and the last
inequality follows from \eqref{full_eq_endpoint_moment_condition}.
Lemma \ref{full_lem_survival_rates}(i) gives
\[
 \lim_{t\to\infty}-\frac1t\log\mbb P_\mu(Z_t>0)=\gamma(\kappa_*).
\]
Dividing the preceding bound by $\mbb P_\mu(Z_t>0)$ and using this
limit yields
\begin{equation}\label{full_eq_endpoint_upper_tail}
 \pi_t((\e^{vt},\infty))\longrightarrow0.
\end{equation}

Choose \(z_0>0\) with \(\mu([z_0,\infty))>0\). For
\(0<\delta<z_0\), monotonicity of
\(z\mapsto\mbb P_z(Z_t>0)\) gives the uniform bound
\[
 \pi_t((0,\delta))
 \le\frac{\mu((0,\delta))}{\mu([z_0,\infty))}
 =:\epsilon_\delta\longrightarrow0\qquad(\delta\downarrow0).
\]
Combining Lemma \ref{full_lem_initial_mass_monotonicity} with
\eqref{full_eq_endpoint_upper_tail}, we obtain
\[
 (1-\epsilon_\delta)R_t(\delta,\lambda)
 \le\int_0^\infty R_t(z,\lambda)\pi_t(\d z)
 \le R_t(\e^{vt},\lambda)+o(1).
\]
The Yaglom limit \cite[Theorem~1.1]{LZZ}, started from $\delta$, gives
\(R_t(\delta,\lambda)\to \Phi_{\kappa_*}(\lambda)\). Lemma
\ref{full_lem_moving_initial_mass} gives
\(R_t(\e^{vt},\lambda)\to
\Phi_{(-{\bf m}-v)/\sigma^2}(\lambda)\).
Let \(\delta\downarrow0\), then
\(v\downarrow-{\bf m}-\kappa_*\sigma^2\), and use
continuity of \(\kappa\mapsto\Phi_\kappa(\lambda)\) at
\(\kappa_*\). The integral converges to
\(\Phi_{\kappa_*}(\lambda)\) for every \(\lambda>0\). The continuity
theorem for Laplace transforms proves the assertion.
\qed}

The preceding results give the following criterion for conditional
convergence.

{
\begin{cor1}
\label{full_cor_survival_criterion}
The conditional laws $\mu_t$ converge weakly to a probability measure on
$(0,\infty)$ if and only if the limit
\begin{equation}\label{full_eq_survival_criterion_limit}
    \lim_{t\to\infty}-\frac1t\log \mbb P_\mu(Z_t>0)
\end{equation}
exists and belongs to $(0,\infty)$. In that case, there is a unique
$\kappa\in(0,\kappa_*]$ such that
\[
    \lim_{t\to\infty}-\frac1t\log \mbb P_\mu(Z_t>0)=\gamma(\kappa),
    \qquad
    \mu_t\Longrightarrow\nu_\kappa.
\]
\end{cor1}

\proof
If $\mu_t$ has a weak limit on $(0,\infty)$, Lemma
\ref{lem_conditional_limit_qsd}(i) and Theorem
\ref{thm_complete_qsd_classification}(i) identify it as $\nu_\kappa$
for some $0<\kappa\leq\kappa_*$.
Lemma \ref{full_lem_conditional_survival_rate} gives the survival rate
$\gamma(\kappa)>0$.

Conversely, suppose that the limit in \eqref{full_eq_survival_criterion_limit}
exists and is strictly positive. By \eqref{full_eq_universal_survival_rate},
this limit belongs to $(0,\gamma(\kappa_*)]$.
Since $\kappa_*\leq-{\bf m}/\sigma^2$,
\[
 \gamma'(r)=-{\bf m}-\sigma^2r>0,\qquad 0\leq r<\kappa_*,
\]
so $\gamma$ is strictly increasing on $[0,\kappa_*]$.
Since $\gamma(0)=0$ and $\gamma$ is continuous, there is a unique
$\kappa\in(0,\kappa_*]$ such that
\[
 \gamma(\kappa)=\lim_{t\to\infty}-\frac1t\log\mbb P_\mu(Z_t>0).
\]
{If $\kappa<\kappa_*$, Proposition
\ref{full_prop_interior_concentration} gives convergence to $\nu_\kappa$.
If $\kappa=\kappa_*$, Lemma \ref{full_lem_survival_rates}(ii) gives
the moment condition \eqref{full_eq_endpoint_moment_condition}, and Proposition
\ref{full_prop_endpoint_attraction} completes the proof.}
\qed}

\subsection{Tail criteria and proof of Theorem \ref{thm_complete_domains_attraction}}\label{sec_tail_criteria}

{It remains to express the survival-rate criterion of Corollary
\ref{full_cor_survival_criterion} in terms of the initial law. We prove the
equivalence of the survival-rate and logarithmic tail conditions for interior
QSDs in Lemma \ref{full_lem_tail_survival_equivalence}. Combining this
equivalence with the endpoint moment criterion from Subsection
\ref{sec_general_initial_laws} then proves Theorem
\ref{thm_complete_domains_attraction}.}

\bglemma\label{full_lem_tail_survival_equivalence}
{Let $0<\kappa<\kappa_*$. Then}
\eqref{full_eq_interior_survival_rate} holds if and only if
\begin{equation}\label{full_eq_tail_limit}
 \lim_{x\to\infty}\frac{-\log\mu([x,\infty))}{\log x}=\kappa.
\end{equation}
\edlemma
\proof
{Suppose first that \eqref{full_eq_tail_limit} holds. Fix
${0<\theta<\kappa}$. Lemma \ref{full_lem_tail_moment} gives
$\int_0^\infty z^\theta\mu(\d z)<\infty$, so
\eqref{full_eq_uniform_survival_moment_bound} yields
\begin{equation}\label{full_eq_tail_survival_upper_bound}
 \mbb P_\mu(Z_t>0)
 \leq C_\theta\left(\int_0^\infty z^\theta\mu(\d z)\right)
       \e^{-\gamma(\theta)t}.
\end{equation}
{For the reverse bound, \eqref{full_eq_tail_limit}, monotonicity, and
Lemma \ref{full_lem_moving_initial_mass} give}
\begin{equation}\label{full_eq_tail_survival_lower_bound}
\begin{aligned}
 \mbb P_\mu(Z_t>0)
 &\geq{\mu([\e^{(-{\bf m}-\kappa\sigma^2)t},\infty))}
       \mbb P_{\e^{(-{\bf m}-{\kappa}\sigma^2)t}}(Z_t>0)\\
 &=\exp\!\left\{-\bigl[{\kappa}(-{\bf m}-{\kappa}\sigma^2)
               +J(-{\bf m}-{\kappa}\sigma^2)\bigr]t+o(t)\right\}\\
 &=\e^{-\gamma({\kappa})t+o(t)}.
\end{aligned}
\end{equation}
Taking logarithms in \eqref{full_eq_tail_survival_upper_bound} and
\eqref{full_eq_tail_survival_lower_bound}, {then letting
$\theta\uparrow\kappa$}, proves \eqref{full_eq_interior_survival_rate}.}

{Conversely, suppose that
\eqref{full_eq_interior_survival_rate} holds.
Lemma \ref{full_lem_survival_rates}(ii) gives finite moments of every
positive order below $\kappa$. Markov's inequality therefore implies
\begin{equation}\label{full_eq_interior_lower_tail}
 \liminf_{x\to\infty}\frac{-\log\mu([x,\infty))}{\log x}
 \geq\kappa.
\end{equation}}
Proposition \ref{full_prop_interior_concentration} and monotonicity of
\(z\mapsto\mbb P_z(Z_t>0)\) give, for every sufficiently small fixed
\(\varepsilon>0\),
\[
 (1-o(1))\mbb P_\mu(Z_t>0)
 \le\mu([\e^{((-{\bf m}-\kappa\sigma^2)-\varepsilon)t},\infty))
                  \mbb P_{\e^{((-{\bf m}-\kappa\sigma^2)+\varepsilon)t}}(Z_t>0).
\]
Indeed the integral over the closed interval in
\eqref{full_eq_interior_concentration} is bounded by the right side.
Taking logarithms, using Lemma \ref{full_lem_moving_initial_mass}, and
setting \(x=\e^{((-{\bf m}-\kappa\sigma^2)-\varepsilon)t}\), we obtain
\[
 \limsup_{x\to\infty}\frac{-\log\mu([x,\infty))}{\log x}
 \le\frac{\gamma(\kappa)-J((-{\bf m}-\kappa\sigma^2)+\varepsilon)}
              {(-{\bf m}-\kappa\sigma^2)-\varepsilon}.
\]
Let \(\varepsilon\downarrow0\). Since
\(\gamma(\kappa)-J((-{\bf m}-\kappa\sigma^2))=\kappa (-{\bf m}-\kappa\sigma^2)\), the upper bound
tends to \(\kappa\). {Together with the lower bound
\eqref{full_eq_interior_lower_tail}, this proves \eqref{full_eq_tail_limit}.}
\qed

\noindent{\it Proof of Theorem \ref{thm_complete_domains_attraction}.~}
{For $0<\kappa<\kappa_*$, Corollary
\ref{full_cor_survival_criterion} and Lemma
\ref{full_lem_tail_survival_equivalence} give
\[
\begin{aligned}
 \mu_t\Longrightarrow\nu_\kappa
 &\quad\Longleftrightarrow\quad
 \lim_{t\to\infty}-\frac1t\log\mbb P_\mu(Z_t>0)=\gamma(\kappa)\\
 &\quad\Longleftrightarrow\quad
 \lim_{x\to\infty}\frac{-\log\mu([x,\infty))}{\log x}=\kappa.
\end{aligned}
\]
For the Yaglom distribution, if $\mu_t\Longrightarrow\nu_{\kappa_*}$,
Lemmas \ref{full_lem_conditional_survival_rate} and
\ref{full_lem_survival_rates}(ii) give the moment condition
\eqref{full_eq_endpoint_moment_condition}. Conversely, that condition
implies attraction to $\nu_{\kappa_*}$ by {Proposition}
\ref{full_prop_endpoint_attraction}.}
Finally, Theorem \ref{thm_complete_qsd_classification} exhausts all QSDs,
so these criteria cover every possible QSD limit.
\qed

{If the logarithmic tail exponent is zero, {then, for every
$0<a<\kappa_*$, we have $\mu([x,\infty))\geq x^{-a}$ for all sufficiently
large $x$. Monotonicity and Lemma \ref{full_lem_moving_initial_mass} give}
\[
 \mbb P_\mu(Z_t>0)
 \geq\e^{-a(-{\bf m}-a\sigma^2)t}
      \mbb P_{\e^{(-{\bf m}-a\sigma^2)t}}(Z_t>0)
 =\e^{-\gamma(a)t+o(t)}.
\]
{Together with} $\mbb P_\mu(Z_t>0)\leq1$, {letting} $a\downarrow0$ {gives}
\[
 \lim_{t\to\infty}-\frac1t\log\mbb P_\mu(Z_t>0)=0.
\]}

\bgremark[Collapse regime]\label{rem_collapse_regime}
In the weakly subcritical case, Theorem
\ref{thm_complete_domains_attraction} shows that attraction to the Yaglom
distribution is equivalent to finiteness of every positive moment of order
below $\kappa_c$.  No logarithmic tail limit is required, and the moment of
order $\kappa_c$ need not be finite.
\edremark
}

\end{document}